\documentclass{article}

\usepackage{epsfig}
\usepackage{times} 
\usepackage[utf8]{inputenc}
\usepackage[T1]{fontenc}
\usepackage{amsmath,stackengine}
\usepackage{amsfonts}
\usepackage{amssymb}
\usepackage{listings}
\usepackage{hyperref}
\usepackage{graphicx}
\usepackage{xcolor}
\usepackage{enumerate}
\usepackage{mathtools}
\usepackage{mathrsfs}
\usepackage{subcaption}
\usepackage{xurl}

\usepackage[linesnumbered,ruled,vlined]{algorithm2e}

\newtheorem{lemma}{Lemma}
\newtheorem{definition}{Definition}
\title{\LARGE \bf
	Output Feedback Reinforcement Learning with Parameter Optimisation for Temperature Control in a Material Extrusion Additive Manufacturing system
}

\author{Eleni Zavrakli*$^{1,2,3}$ , Andrew Parnell $^{1,2,3}$ and Subhrakanti Dey$^{4}$%
	\thanks{*Corresponding author: \tt\small eleni.zavrakli@mu.ie}%
	\thanks{$^{1}$Department of Mathematics and Statistics,
		Maynooth University, Co. Kildare, Ireland}%
	\thanks{$^{2}$I-Form Advanced Manufacturing Research Centre, Ireland}%
	\thanks{$^{3}$Hamilton Institute, Maynooth University, Co. Kildare, Irealand}%
	\thanks{$^{4}$Division of Signals and Systems, Department of Electrical Engineering, Uppsala University, Sweden}
}

\begin{document}
	
	\maketitle
	\begin{abstract}
		With the rapid development of Additive Manufacturing (AM) comes an urgent need for advanced monitoring and control of the process. Many aspects of the AM process play a significant role in the efficiency, accuracy and repeatability of the process, with temperature regulation being one of the most important ones. In this work, we solve the problem of optimal tracking control for a state space temperature model of a Big Area Additive Manufacturing (BAAM) system. In particular, we address the problem of designing a Linear Quadratic Tracking (LQT) controller when access to the exact system state is not possible, except in the form of measurements. We initially solve the problem with a model-based approach based on reinforcement learning concepts, with state estimation through an observer. We then design a model-free reinforcement-learning based controller with an internal state estimation step and demonstrate its performance through a simulator of the systems' behaviour. Our results showcase the possibility of achieving comparable results while learning optimal policies directly from process data, without the need for an accurate, intricate model of the process. We consider this outcome to be a significant stride towards autonomous intelligent manufacturing.
	\end{abstract}
	

\maketitle
\textbf{Acknowledgements}:
This publication has emanated from research supported by a research grant from Science Foundation Ireland (SFI) under Grant Number 16/RC/3872. For the purpose of Open Access, the author has applied a CC BY public copyright licence to any Author Accepted Manuscript version arising from this submission.

\textbf{Competing interests}: The authors have no competing interests to declare that are relevant to the content of this article.

\textbf{Code availability}: The code used for our simulations was written in R and can be found on: 
\url{https://github.com/elenizavrakli/Output-feedback-RL-temperature-control}

\section{Introduction}
Additive Manufacturing (AM) is becoming one of the most popular manufacturing techniques, with its applications ranging from aerospace to bio-medicine to household items \cite{BERMAN2012155,BUCHANAN2019332,DILBEROGLU2017545}. Its appeal lies in the ability to build very complicated structures without the need for a mold and without the waste of a lot of materials. There is an array of different AM techniques, using different materials, and targeting different problems. For example, in the industrial scale, metal-based AM is being adopted at an increasing scale due to its cost-effective and energy-saving nature \cite{frazier2014metal,uriondo2015present,gu2012laser}. Material extrusion (MEX) on the other hand, a technique using different polymers such as PLA, allows consumers to design and manufacture their own prototypes at home, making it more appealing to hobbyists and small business owners, neither of which would normally have access to traditional manufacturing techniques and equipment. As AM becomes increasingly prevalent, the need naturally arises for the process to be optimised in terms of safety, repeatability and cost effectiveness. This objective can be addressed by efficient closed-loop control algorithms, which are currently lacking in the AM space \cite{mercado2020additive}. With many different aspects of AM processes needing to be monitored and controlled, developing effective feedback controllers to target them is a very active field of study \cite{lhachemi2019augmented}. 

Control Theory \cite{anderson2007optimal,ogata2010modern} offers a wide variety of methods for designing efficient controllers, normally when a mathematical model describing the process is available. In the AM context, proportional integral derivative (PID) controllers have been widely used to control various aspects of the process \cite{previdi2006design,farshidianfar2016real-time,HU200351,kruth2007feedback,8301604}. In AM, similarly to many other physical systems, obtaining accurate models of the process is very challenging, which is an active research field on its own \cite{bikas2016additive}. With that in mind, the development of controllers tuned through information learned from process data is being increasingly investigated. Reinforcement Learning (RL) \cite{si2004handbook,sutton2018reinforcement}
commonly refers to a category of data-centric Machine Learning techniques that learn optimal strategies within intricate and unpredictable settings by maximizing rewards (or minimizing costs). In the field of Control Theory, decision-making problems expressed as optimising strategies are formulated using Dynamic Programming \cite{bellman1966dynamic,bertsekas2012dynamic,bertsekas1996neuro}, which shares several commonalities with RL but employs distinct terminology. While certain RL methods have been designed deterministically, assuming a priori knowledge of the system dynamics, the majority of RL applications are rooted in data-driven, model-free approaches. In recent years, as RL has gained popularity, researchers have employed it to address various issues within AM \cite{8362941,wasmer2019situ,alicastro2021reinforcement,patrick2018reinforcement,dharmawan2020model,ogoke2021thermal}. However, none of these papers use RL as a tool to design optimal feedback controllers. For more references we refer to our previous work \cite{zavrakli2023data}, where we investigated the problem of temperature control in an AM system, when the exact system state is directly available. The problem was solved using both standard Control Theory approaches and RL techniques, both in a model-based and a data-driven manner.  

In this paper, we address the problem of temperature control for a MEX AM system without access to the exact system state, but with available sensor measurements only. Specifically, we control the temperatures within the extruder of a Big Area Additive Manufacturing (BAAM) system. We design a controller that solves the Linear Quadratic Tracking (LQT) problem, derived directly from the process model. However, we ultimately aim to have a controller that can learn optimal policies directly from process data. We utilise Reinforcement Learning techniques to design a data-driven LQT controller \cite{kiumarsi2014reinforcement}.  In the model-based case, the state can be estimated with the help of an observer \cite{luenberger1971introduction}. In the case of the data-driven controller, we incorporate a data-driven state-estimation step in the controller, using a record of past inputs and outputs from the system only, inspired by the work done in \cite{9146286}. Furthermore, we optimise the parameters used to define the LQT problem using Bayesian Optimisation \cite{frazier2018tutorial,shahriari2015taking}. As there is no systematic way to choose these parameters, Bayesian Optimisation has been explored in the literature due to its ability to handle optimising complex objectives when trial-and-error approaches do not suffice \cite{calandra2016bayesian,marco2016automatic,miyamoto2018automatic}. Finally, while setting up our simulations, we provide detailed guidelines on how to generate appropriate data to train the data-driven controller. These data generation tools can be useful in many settings, in particular when one needs to  train a learning algorithm for dynamical systems.

The structure of our paper is as follows: In Section \ref{section_sys_intro}, we introduce the BAAM system we wish to control and define the optimisation problem to be solved. Section \ref{section_model_based} deals with the design of the model-based state observer and controller, while in Section \ref{section_data_driven} we present the design of the RL-based data-driven controller which includes state estimation through input-output data. Section \ref{section_BO} introduces the main concepts of Bayesian Optimisation and Section \ref{section_results} presents the results of our simulations, using both controllers, before and after Bayesian Optimisation. Finally, in Section \ref{section_conclusion} we provide a thorough discussion of our findings.

\section{Linear Quadratic Output Feedback Control for a state-space system}\label{section_sys_intro}

Faulty prototypes often occur in AM processes due to issues related to temperature regulation in various subprocesses in AM. More specifically, the cause of these faults may be the insufficient or excessive melting of the material, the unsuccessful binding between layers or the non-uniform heating of the build-plate which can stop the print from adhering properly. 

\begin{figure}[ht]
    \centering
    \includegraphics[width=\textwidth]{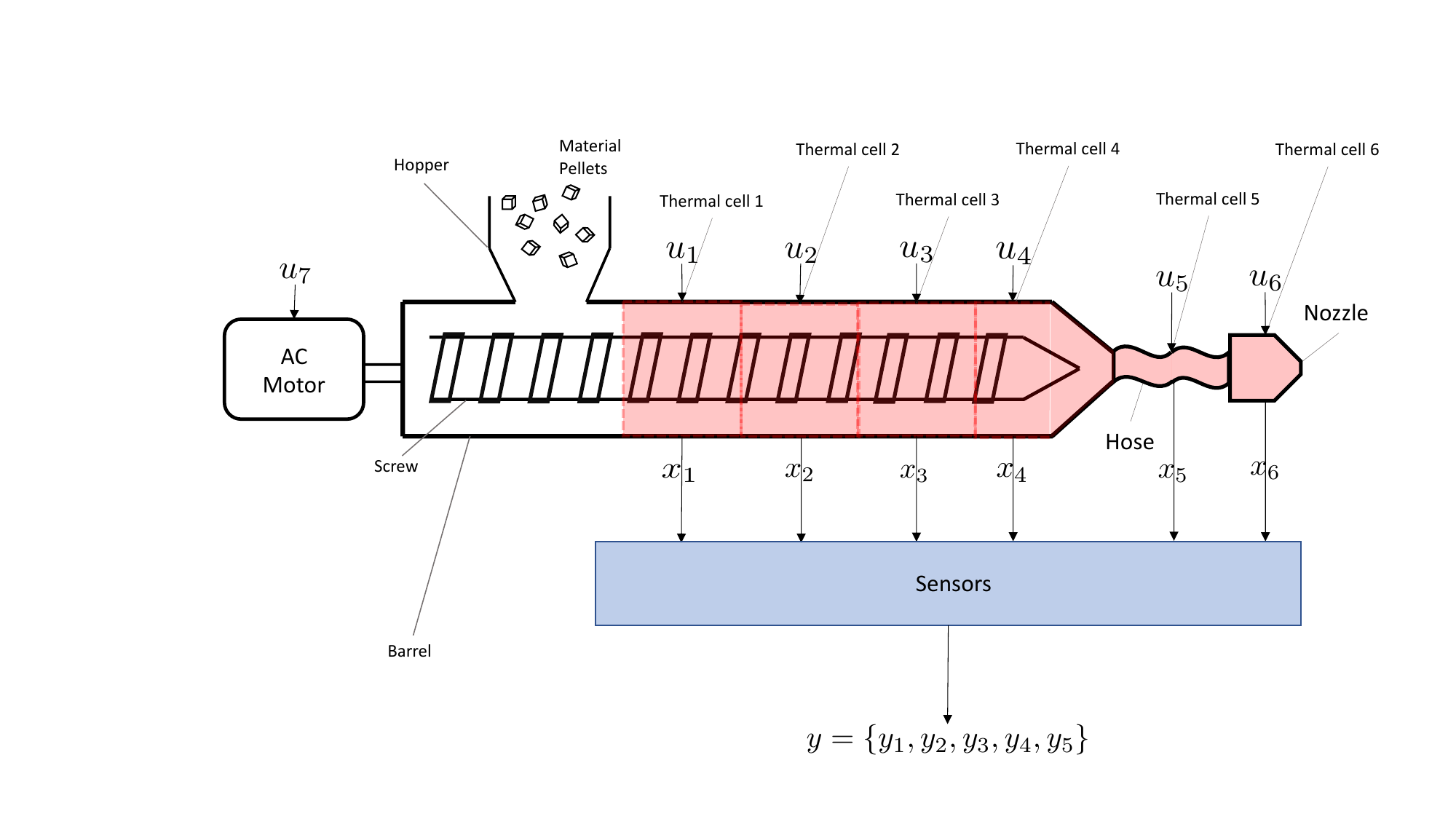}
    \caption{Heating system for a BAAM extruder. The inputs include the heat applied in each thermal cell and the input provided by the AC motor rotating the screw. The system states represent the temperatures within the extruder, which are measured by the sensors.}
    \label{schematic}
\end{figure}

In this work, we study the heating system of a Material Extrusion (MEX) system. In particular, we are focusing on the temperatures within the extruder of a Big Area Additive Manufacturing (BAAM) system. The extruder consists of 5 main parts, the hopper, the barrel, the hose, the nozzle and the AC motor. The hopper is where the material is inserted and gets fed to the barrel containing the screw which gets rotated by the motor. The screw pushes the material towards the hose and finally to the nozzle. There are four heaters placed along the barrel, one heater in the hose and one in the nozzle. The reason behind the use of multiple heaters in BAAM systems is to address the large volume of material that gets processed by the system. The areas where the heat is applied can be viewed as thermal cells, each influencing their neighbouring cells and all getting influenced by the motor. We adopt a discrete time, linear state-space representation of the system inspired by the system model obtained in \cite{gootjes2017applying} through system identification methods using real process data. We make the further assumption that the temperatures within the extruder are not directly available, instead a set of measurements can be obtained through a number of sensor measurements of the process.
\begin{eqnarray}\label{original_sys}
    x(t+1)&=&Ax(t)+Bu(t) \nonumber\\ 
    y(t) &=& Cx(t)
\end{eqnarray}
where $t\geq t_0$ is the discrete time step, $x(t)\in \mathbb{R}^{n}$ is the system state at time $t$, $u(t)\in \mathbb{R}^{m}$ is the system input at $t$ and $y(t)\in \mathbb{R}^{p}$ is the system output obtained as measurements. Matrices $A\in \mathbb{R}^{n\times n}, B\in \mathbb{R}^{n \times m}$ and $C\in \mathbb{R}^{p\times n}$ are the state, input and measurement matrices respectively. For our system, there are $m=7$ inputs, the first 6 representing the load applied to each thermal cell by the corresponding heater and the last relaying the input applied by the motor. The system state is a vector of length $n=6$ representing the temperature within each thermal cell. We assume that the sensor system obtains $p=5$ measurements of the internal states of the system. A visual representation of the system setup can be seen in Figure \ref{schematic}.

The temperature control problem is defined as the infinite horizon optimal tracking problem, where the goal is to obtain the optimal input (or control) $u^*$ such that the system output accurately tracks the reference signal $r(t)\in \mathbb{R}^p$, generated according to matrix $F\in \mathbb{R}^{p\times p}$
\begin{equation}\label{reference_sys}
    r(t+1)=Fr(t).
\end{equation}
In most applications $r(t)$ is chosen to be constant, namely $F=I_p$ or $r(t)=r\in \mathbb{R}$ for all $t$. However, the general, time-varying reference signal \eqref{reference_sys} can also be adopted in our formulation.
The control objective can be formally expressed as the minimisation of the performance index defined as the total sum of discounted costs over an infinite horizon, namely:
\begin{eqnarray}\label{cost_func}
   \min_\mathbf{u} J(x(0),\mathbf{u})&=  &\min_\mathbf{u} \sum_{t=0}^{\infty} \gamma^{t} \left\lbrace \lbrack y(t)-r(t)\rbrack^T Q\lbrack y(t)-r(t)\rbrack + u^T(t)Ru(t)  \right\rbrace \nonumber\\ 
    &= &\min_\mathbf{u} \sum_{t=0}^{\infty} \gamma^{t} \left\lbrace \lbrack Cx(t)-r(t)\rbrack^T Q\lbrack Cx(t)-r(t)\rbrack+ u^T(t)Ru(t)  \right\rbrace \nonumber\\
\end{eqnarray}
where $0<\gamma \leq 1$ is the discount factor, putting greater emphasis on short term costs. $Q\in \mathbb{R}^{p\times p}$ and $R\in \mathbb{R}^{m\times m}$ are positive definite symmetric matrices applying weight to the tracking error and the input cost respectively. $x(0)$ gives rise to all future states $x(t)$ and we define $\mathbf{u}$ as the infinite sequence containing all inputs $\{ u(0), u(1), u(2),\dots \}$, also known as a \textit{policy}. This means that the minimisation problem is equivalent to the search for an optimal policy $\mathbf{u}^*$. In the Reinforcement Learning literature \cite{bellman1966dynamic,bertsekas2012dynamic,sutton2018reinforcement}, the minimisation problem is addressed through the principle of optimality:

	\begin{definition}[Principle of optimality] \cite{bellman1966dynamic}
		An optimal policy has the property that whatever the initial state and initial input are, the remaining inputs must constitute an optimal policy with regard to the state resulting from the first decision.
	\end{definition}
To utilise the principle of optimality, we define the \textit{Value Function} as the cost-to-go at time $t$, namely the sum of discounted future costs from time $t$ onwards.
\begin{eqnarray}\label{perf_index}
    V(x(t),\mathbf{u}_t)&=& \sum_{k=t}^{\infty} \gamma^{k-t} \left\lbrace \lbrack y(k)-r(k)\rbrack^T Q\lbrack y(k)-r(k)\rbrack + u^T(k)Ru(k)  \right\rbrace \nonumber\\ 
    &=& \sum_{k=t}^{\infty} \gamma^{k-t} \left\lbrace \lbrack Cx(k)-r(k)\rbrack^T Q\lbrack Cx(k)-r(k)\rbrack + u^T(k)Ru(k)  \right\rbrace
\end{eqnarray}
where $\mathbf{u}_t=\{ u(t), u(t+1), u(t+2),\dots \}$. The Value function allows for the principle of optimality to be expressed mathematically as the \textbf{Bellman optimality equation} \cite{bertsekas2012dynamic, sutton2018reinforcement}. Let $c(x(t),u(t))=\lbrack Cx(t)-r(t)\rbrack^T Q\lbrack Cx(t)-r(t)\rbrack + u^T(t)Ru(t) $ be the cost at time $t$. Then for the optimal policy $\mathbf{u}^*$, the principle of optimality yields
\begin{equation}\label{bellman_general}
     V(x(t),\mathbf{u}^*_t)= c(x(t),u^*(t))+\gamma  V(x(t+1),\mathbf{u}^*_{t+1})
\end{equation}
for all $t$.
In the following sections we explore two different approaches to solving the optimal tracking problem. The first is based on the design of a state observer while the second utilises a record of past inputs and outputs to estimate the internal state. While the observer design requires knowledge of the system dynamics, the second approach can be adapted to be completely model free and data-driven.

\section{State Observer based Output Feedback}\label{section_model_based}

An important issue in the design of efficient controllers is the lack of access to the exact system state $x(t)$. Assuming knowledge of the system dynamics \eqref{original_sys}, a state observer can be implemented to provide an informed estimate of the state, that can then be used for control design. This observer is also known as a closed loop observer or a Luenberger observer \cite{luenberger1971introduction,ogata2010modern}. 

In practice, the Luenberger observer is another state space system, made up of the same state, input and output matrices as the original system, where the state equation has an additional term associated with the error between the system output and the observer output. 
\begin{eqnarray}\label{observer_equations}
\hat{x}(t+1)&=&A\hat{x}(t)+Bu(t) +L\left\lbrack y(t)-\hat{y}(t)\right\rbrack \nonumber \\ 
\hat{y}(t) &=& C\hat{x}(t)
\end{eqnarray}
The matrix $L \in \mathbb{R}^{n \times p}$ is known as the observer gain. The goal is to have an observer that emulates the behaviour of the original system as closely as possible. In control theory terms, we need to find a gain matrix $L$ such that the matrix $A-LC$ is stable \cite{ogata2010modern}. This problem can be approached as a pole placement problem, where the objective is to place the poles of $A-LC$ inside the unit circle. For that to be possible, the system needs to be  \textit{completely state observable} \cite{ogata2010modern}. The observability condition requires for the column rank of the matrix
\begin{equation}\label{observability}
    \mathcal{O} = 
    \begin{bmatrix}
    C \\ CA \\ CA^2 \\  \vdots \\ CA^{n-1}
    \end{bmatrix}
\end{equation}
to be $n$. If this condition is satisfied, the system state can be fully determined from the system inputs and outputs. One approach to solve the pole placement problem is to choose a small positive number $\tau$ and calculate the observer gain through:
\begin{equation}
    L=AC^T\Psi^{-1} \text{ where } \Psi=CC^T+\tau I_p
\end{equation}
Once $L$ is determined, the state estimate $\hat{x}$ can be used in place of $x$ to design a controller. This means that the value function \eqref{perf_index} can be rewritten as follows: 
\begin{equation}\label{perf_index2}
    V_{obs}(x(t),\mathbf{u}_t)= \sum_{k=t}^{\infty} \gamma^{k-t} \left\lbrace \lbrack C\hat{x}(k)-r(k)\rbrack^T Q\lbrack C\hat{x}(k)-r(k)\rbrack + u^T(k)Ru(k) \right\rbrace
\end{equation}

For the optimal control problem to be successfully solved, the system needs to be \textit{controllable} \cite{ogata2010modern}. System \eqref{original_sys} is said to be controllable when the matrix
\begin{equation}\label{controllability}
    \mathcal{C} = \begin{bmatrix}B& AB & A^2B & \cdots & A^{n-1}B\end{bmatrix}
\end{equation}
known as the controllability matrix, has full row rank $n$. The controllability condition implies that the system can be controlled to have the desired behaviour within a finite number of steps, through an appropriate choice of control function $u$.

For the infinite horizon LQT problem \eqref{perf_index2}, a controller can be designed by adopting the approach in \cite{kiumarsi2014reinforcement}. Consider the augmented state 
\begin{equation}\label{aug_state}
    X_{obs}(t)=\begin{bmatrix} \hat{x}(t) \\ r(t)
    \end{bmatrix}
\end{equation}
and the corresponding system dynamics
\begin{equation}\label{aug_sys_dynamics}
    X_{obs}(t+1)= TX_{obs}(t)+B_1u(t).
\end{equation}
where $T=\begin{bmatrix} A& 0 \\ 0&F\end{bmatrix}$ and $B_1=\begin{bmatrix} B \\ 0\end{bmatrix}$. Using \eqref{aug_state}, \eqref{perf_index2} can be rewritten as
\begin{equation}\label{perf_index3}
    V_{obs}(X(t),\mathbf{u}_t)= \sum_{k=t}^{\infty} \gamma^{k-t}\lbrace X_{obs}^T(k) Q_1X_{obs}(k) + u^T(k)Ru(k)  \rbrace
\end{equation}
where 
\begin{equation}\label{aug_weighting_q}
Q_1=  \begin{bmatrix}C^TQC & -C^TQ \\ -QC & Q \end{bmatrix}.
\end{equation}
The optimal policy is in feedback form of the augmented system state
\begin{equation}\label{policy_state}
    u^*(t)=-K X_{obs}(t).
\end{equation}
where $K$ is the optimal gain. $K$ is commonly obtained through the solution of an Algebraic Riccati Equation (ARE) which can be a challenging task and computationally expensive. An alternative approach to determining $K$, particularly useful for large scale models, is through Reinforcement Learning (RL) techniques, namely, by utilising the Bellman equation \cite{bertsekas1996neuro,sutton2018reinforcement,si2004handbook}. In order to form the Bellman equation, the following lemma is needed.
\begin{lemma}\cite{kiumarsi2014reinforcement} \label{lemma}
For the optimal tracking problem with value function of the form \eqref{perf_index2} and reference generated by $r_{obs}(t+1)=F_{obs}x_{obs}(t)$, then for any stabilising policy of the form \eqref{policy_state}, the value function can be written in quadratic form as 
\begin{equation}\label{quad_form_value_function_obs}
    V_{obs}(x(t),\mathbf{u}_t)= V_{obs}(X(t))=\frac{1}{2}X_{obs}^T(t)PX(_{obs}t)
\end{equation}
for some matrix $P=P^T>0$.
\end{lemma}
Utilising the quadratic form of the value function, it can be rewritten as the Bellman equation 
\begin{equation}\label{bellman}
    X_{obs}^T(t)PX_{obs}(t)=X_{obs}^T(t)Q_1X_{obs}(t)+u^T(t)Ru(t) +\gamma X_{obs}^T(t+1)PX_{obs}(t+1).
\end{equation}
Consider a stabilising policy $\hat{K}^*$ in \eqref{bellman}. A matrix $K$ is said to be stabilising when the matrix $A-BK$ is stable, namely has eigenvalues with strictly negative real parts \cite{anderson2007optimal}. Using the augmented system dynamics \eqref{aug_sys_dynamics}, we obtain the Lyapunov equation
\begin{equation}\label{lyapunov}
    P=Q_1 +(\hat{K}^*)^TR\hat{K}^*+\gamma (T+B_1\hat{K}^*)^TP(T+B_1\hat{K}^*).
\end{equation}
Iterative solutions of the Lyapunov equation while updating the optimal control estimate $\hat{u}^*=-\hat{K^*}X_{obs}$ can effectively lead to the solution to the Optimal Tracking Problem. A step by step description of the RL based solution to the LQT problem using state estimation through a state observer is given in Algorithm \ref{observer_algorithm}.

\begin{algorithm}
		\SetAlgoLined
		\KwResult{State estimate $\hat{x}(t)$, Optimal control gain $K^*$, optimal controller $u^*(t)$}
		\KwIn{State input and output matrices $A,B,C$, Weighting matrices for performance index $Q,R$, Initial state estimate $\hat{x}(t_0)$, small positive number $\tau$, Initial reference signal $r(t_0)$,  Discount factor $\gamma$, Initial control policy $K^0$, Error threshold $\epsilon$}
		\begin{enumerate}
		    \item Create augmented system matrices $T=\begin{bmatrix}
	A & 0 \\ 0 & F
	\end{bmatrix},B_1=\begin{bmatrix}
	B \\ 0 
	\end{bmatrix}$.
	\item Augment weighting matrix $Q_1=  \begin{bmatrix}C^TQC & -C^TQ \\ -QC & Q \end{bmatrix}.$
	\item Determine optimal control gain $K^*$:
	\begin{enumerate}[(i)]
	    \item Set $j=0$.
	    \item Solve the LQT Lyapunov equation $P^{j+1}=Q_1+(K^j)^T RK^j +\gamma (T-B_1K^j)^TP^{j+1}(T-B_1K^j)$.
	\item Compute the control gain estimate 
		$K^{j+1}= (R+\gamma B_1^T P^{j+1} B_1)^{-1}$  $\gamma B_1^TP^{j+1}T$.
		\item Set $j=j+1$.
	\end{enumerate}
	Repeat steps $(ii)-(iv)$ while $\vert K^{j}-K^{j-1}\vert>\epsilon$. 
	\item Set $K^*=K^j$.
	\item Determine and apply optimal controller to system while estimating the state through the observer:
		\begin{enumerate}[(a)]
            \item Set $t=t_0$.
            \item Design augmented state $X_{obs}(t)=\begin{bmatrix} \hat{x}(t) \\ r(t)\end{bmatrix}$
		    \item Obtain system input $u^*(t)=-K^*X_{obs}(t)$ 
		    \item Obtain output estimation $\hat{y}(t)=C\hat{x}(t)$.
		    \item Observe system output $y(t)=Cx(t)$.
		    \item Determine $L=AC^T \Psi^{-1}$ where $\Psi = CC^T + \tau I_p$.
		    \item State estimate $\hat{x}(t+1)=A\hat{x}(t)+Bu^*(t) +L(y(t)-\hat{y}(t))$.
		    \item System state update (not measurable) $x(t+1)=Ax(t)+Bu^*(t)$
            \item Reference signal update $r(t+1)=Fr(t)$
		    \item Set $t=t+1$.
		\end{enumerate}
		Repeat steps $(b)-(j)$ for the desired length of the experiment or simulation.
		\end{enumerate}
		\caption{ \label{observer_algorithm} LQT control with state estimation through Luenberger observer}
	\end{algorithm}   

    If the initial conditions given to the observer are the same or close in value to the true initial conditions of the system, the observer estimates the state very closely within a few time steps. The further away the initial values are, the longer it takes for the observer to estimate the state correctly, so the longer it takes for the system's behaviour to be optimised, as the controller is designed based on the state estimates. 

    While the observer is a powerful estimator of the state, it includes an important assumption which is access to the system dynamics, namely matrices A, B and C. This is an unrealistic assumption in many applications, where a model of the process is either not available or not accurate enough. In the next section, we explore another approach for control design with state estimation. 

    \section{Output Feedback using Input-Output data}\label{section_data_driven}
    Consider the performance index \eqref{perf_index} to be minimised with $r$ generated as in \eqref{reference_sys}. Similarly to the previous section, we augment the state with the reference and construct the augmented system dynamics:
    \begin{equation} \label{aug_sys_i_o}
    X(t+1)= \begin{bmatrix} x(t+1) \\ r(t+1) 
    \end{bmatrix}= \begin{bmatrix} A& 0 \\ 0&F
    \end{bmatrix}\begin{bmatrix} x(t) \\ r(t)
    \end{bmatrix}+\begin{bmatrix} B \\ 0
    \end{bmatrix}u(t) = TX(t)+B_1u(t).
\end{equation}
The value function can be written in terms of $X$ as
\begin{equation}\label{perf_index4}
    V(X(t),\mathbf{u}_t)= \sum_{k=t}^{\infty} \gamma^{k-t}\lbrace X^T(k) Q_1X(k) + u^T(k)Ru(k)  \rbrace
\end{equation}
where $Q_1$ defined as in \eqref{aug_weighting_q}.
Assuming a feedback control policy, \ref{lemma} allows for \ref{perf_index4} to be written as
\begin{equation}\label{quad_form_value_function}
    V(x(t),\mathbf{u}_t)= V(X(t))=\frac{1}{2}X^T(t)P_1X(t)
\end{equation}
which gives rise to the Bellman equation 
\begin{equation}\label{bellman2}
    X^T(t)P_1X(t)=X^T(t)Q_1X(t)+u^T(t)Ru(t)+\gamma X^T(t+1)P_1X(t+1).
\end{equation}

Unlike the previous section, in addition to not having access to the internal system state, we no longer assume knowledge of the system dynamics. Instead, we estimate the state through a record of past inputs and outputs, under the assumption that the system \eqref{original_sys} is controllable and observable. Consider a time horizon $N$ and recursively express the state $x(t)$ in terms of past inputs $u(t-1), \dots, u(t-N)$
\begin{eqnarray}\label{recursive state expression}
        x(t) &=& Ax(t-1)+Bu(t-1)\nonumber \\ 
        &=& A^2x(t-1)+ ABu(t-2)+Bu(t-1) \nonumber \\ && \dots \nonumber \\ 
        &=&
        A^N x(t-N) + \begin{bmatrix}
         B & AB &A^2 B & \dots & A^{N-1}B 
    \end{bmatrix} \begin{bmatrix}
        u(t-1) \\ u(t-2) \\ \vdots \\ u(t-N)
    \end{bmatrix} \nonumber\\ &:=&A^N x(t-N)  + U_N \bar{u}(t-1,t-N).
\end{eqnarray}
Similarly, the system output can be expressed recursively as
\begin{equation} \label{recursive output expression}
    y(t)=Cx(t)=CA^Nx(t-N)+C U_N \bar{u}(t-1,t-N).
\end{equation}
Consider a sequence of N outputs
\begin{equation}
 \bar{y}(t-1,t-N)= \begin{bmatrix}
        y(t-1) \\ y(t-2) \\ y(t-3) \\ \vdots \\ y(t-N)
    \end{bmatrix}    
\end{equation}
    which can be expressed in terms of $x(t-N)$ and $\bar{u}(t-1,t-N)$ as
    \begin{eqnarray}\label{bary equation}
        \bar{y}(t-1,t-N) & =&  \begin{bmatrix}
        CA^{N-1} \\  CA^{N-2} \\ \vdots \\ CA \\ C
    \end{bmatrix}
        x(t-N) \nonumber\\ && +\begin{bmatrix}
         0 & CB & CAB  & \dots & CA^{N-2}B \\
          0 & 0 & CB  & \dots & CA^{N-3}B \\
           \vdots & \vdots & \vdots  & \ddots & \vdots \\
           0 & 0 &  0 & \dots & CB \\
           0 & 0& 0& \dots &0
    \end{bmatrix} \begin{bmatrix}
        u(t-1) \\ u(t-2) \\ \vdots \\ u(t-N)
    \end{bmatrix} \nonumber\\ 
 &:=& W_N x(t-N) + D_N \bar{u}(t-1,t-N)
    \end{eqnarray}
    Matrix $W_N$ has full rank for $N\geq K$ where $K$ is the rank of the observability matrix \eqref{observability}. Since we assume our system to be observable, then $K=n$. If $N$ is chosen to be equal to $n$, matrix $U_N$ is equal to the controllability matrix $\mathcal{C}$ \ref{controllability} and $W_N$ is equal to the observability matrix $\mathcal{O}$ \ref{observability} flipped from top to bottom. 
    
    We aim to express the state $x(t)$ in terms of past inputs and outputs. Based on \ref{recursive state expression}, $x(t-N)$ needs to be replaced by an expression in terms of $u(t-1), \dots, u(t-N)$ and $y(t-1), \dots, y(t-N)$. This can be achieved by solving \ref{recursive output expression} for $x(t-N)$. To that end, consider the generalised Moore-Penrose inverse \cite{moore1920reciprocal,penrose_1955} of $W_N$, $W^+_N=(W^T_N W_N)^{-1}W_N^T$. Then 
    \begin{equation}
        x(t-N) =  W^+_N(\bar{y}(t-1,t-N)-D_N\bar{u}(t-1,t-N)). 
    \end{equation}
    Considering that $r(t)=F^Nr(t-N)$, the augmented state $X(t)$ can be written as:
    \begin{eqnarray}\label{aug_state_estimation}
  X(t) &=&
  \begin{bmatrix}
      A^N & 0 \\ 0 & F^N
  \end{bmatrix}  \begin{bmatrix}
      x(t-N) \\ r(t-N)
  \end{bmatrix}  + \begin{bmatrix}
      U_N \\ 0
  \end{bmatrix}\bar{u}(t-1,t-N) \nonumber\\ 
  &=&
  \begin{bmatrix}
      A^N & 0 \\ 0 & F^N
  \end{bmatrix}  \begin{bmatrix}
      W^+_N(\bar{y}(t-1,t-N)-D_N\bar{u}(t-1,t-N)) \\ r(t-N)
  \end{bmatrix}   \nonumber\\  &&+ \begin{bmatrix}
      U_N \\ 0
  \end{bmatrix}\bar{u}(t-1,t-N) \nonumber\\ 
  &=& \begin{bmatrix}
      U_N-A^N W^+_N D_N & A^NW^+_N & 0 \\ 0 & 0 & F^N 
  \end{bmatrix}\begin{bmatrix}
       \bar{u}(t-1,t-N)\\  \bar{y}(t-1,t-N) \\ r(t-N)
  \end{bmatrix}   \nonumber \\
  &:= & M \begin{bmatrix}
       \bar{u}(t-1,t-N)\\  \bar{y}(t-1,t-N) \\ r(t-N)
  \end{bmatrix}
\end{eqnarray}
The above result establishes that the system state can be expressed in terms of past inputs and outputs with the help of matrix $M$. Let $Z(t)=\begin{bmatrix}
    X(t)\\u(t)
\end{bmatrix}$ and consider the quadratic form of the value function \ref{quad_form_value_function}. It can be rewritten as
\begin{equation}\label{quad_form_aug}
    V(Z(t))=\frac{1}{2}Z^T(t)HZ(t)
\end{equation}
with $H=\begin{bmatrix} Q_1+\gamma T^TP_1T & \gamma T^TP_1B_1 \\ \gamma B_1^TP_1T &R+\gamma B_1^TP_1B_1
    \end{bmatrix}$.
Using \ref{aug_state_estimation} in \ref{quad_form_aug} we obtain
  \begin{eqnarray}
           V(Z(t))&=&\frac{1}{2}\begin{bmatrix}
          X(t)\\u(t) 
      \end{bmatrix}^T H \begin{bmatrix}
          X(t)\\u(t) 
      \end{bmatrix} \nonumber \\  &=& \frac{1}{2}\begin{bmatrix}
          M \begin{bmatrix}
       \bar{u}(t-1,t-N)\\  \bar{y}(t-1,t-N) \\ r(t-N)
  \end{bmatrix}  \\u(t) 
      \end{bmatrix}^T H \begin{bmatrix}
          M \begin{bmatrix}
       \bar{u}(t-1,t-N)\\  \bar{y}(t-1,t-N) \\ r(t-N)
  \end{bmatrix}  \\u(t) 
      \end{bmatrix} 
  \end{eqnarray}
 Set $\bar{Z}(t)=\begin{bmatrix}
       \bar{u}(t-1,t-N)\\  \bar{y}(t-1,t-N) \\ r(t-N) \\u(t)
  \end{bmatrix}$ and $\bar{H}= \begin{bmatrix}
      M &0\\0& I_{m}
  \end{bmatrix}^T H \begin{bmatrix}
      M &0\\0& I_{m}
  \end{bmatrix}$. Then we can rewrite \ref{quad_form_aug} as
  \begin{eqnarray} \label{final_quad_form}
     V(\bar{Z}(t))& =& \frac{1}{2}\bar{Z}^T(t) \bar{H} \bar{Z}(t) \nonumber\\  &=& {\small\frac{1}{2} \begin{bmatrix}
      \bar{u}(t-1,t-N+1)\\  \bar{y}(t,t-N+1) \\ r(t-N+1) \\u(t)
  \end{bmatrix}^T 
  \begin{bmatrix}
      H_{\bar{u}\bar{u}} & H_{\bar{u}\bar{y}} &H_{\bar{u}r} & H_{\bar{u}u} \\ 
      H_{\bar{y}\bar{u}} & H_{\bar{y}\bar{y}}&H_{\bar{y}r} & H_{\bar{y}u} \\ 
      H_{r\bar{u}} &  H_{r\bar{y}} & H_{rr} & H_{ru} \\ 
      H_{u\bar{u}} &  H_{u\bar{y}} & H_{ur} &H_{uu}
  \end{bmatrix} 
  \begin{bmatrix}
      \bar{u}(t-1,t-N+1)\\  \bar{y}(t,t-N+1) \\ r(t-N+1) \\u(t)
  \end{bmatrix}} \nonumber \\
  \end{eqnarray}
The decomposition of $\bar{H}$ in blocks of the form $H_{ij}$ is according to the dimensions of each set of $i,j$. For example block $H_{\bar{u}\bar{u}}$ will be a $Nm\times Nm$ matrix and block $H_{u\bar{y}}$ will be a $m\times Np$ matrix.
  
Due to the quadratic form of the value function, the optimality condition $\frac{\partial V}{\partial u}=0$  is sufficient to ensure its minimisation. Solving for $u(t)$ yields
\begin{equation}\label{optimal_control_func}
        u(t) = - H_{uu}^{-1} ( H_{u\bar{u}} \bar{u}(t-1,t-N)+ H_{u\bar{y}}\bar{y}(t-N,t-N) + H_{ur} r(t-N)).
    \end{equation}

   If the system model is available, \ref{optimal_control_func} can be used to directly determine the optimal LQT controller using input-output information. The kernel matrix $\bar{H}$ can be directly determined using 
   \begin{eqnarray}\label{kernel}
       \bar{H} &=\begin{bmatrix}
      M &0\\0& I_{m}
  \end{bmatrix}^T H \begin{bmatrix}
      M &0\\0& I_{m}
  \end{bmatrix}  = \begin{bmatrix}
      U_N-A^N W^+_N D_N & A^NW^+_N & 0 & 0 \\ 0 & 0 & F^N &0 \\ 0&0&0&I_m
  \end{bmatrix}^T \nonumber \\ &\begin{bmatrix} Q_1+\gamma T^TP_1T & \gamma T^TP_1B_1 \\ \gamma B_1^TP_1T &R+\gamma B_1^TP_1B_1
    \end{bmatrix}\begin{bmatrix}
      U_N-A^N W^+_N D_N & A^NW^+_N & 0 & 0 \\ 0 & 0 & F^N &0 \\ 0&0&0&I_m
  \end{bmatrix}\nonumber \\
   \end{eqnarray}
 where $P_1$ is obtained by solving the standard LQT problem as in \ref{lyapunov}.     

We are interested in the case when a reliable model of the system is not available, meaning matrices $A,B$ and $C$ are unknown. In this case, the kernel matrix $\bar{H}$ needs to be estimated strictly from data. To that end, consider the Bellman equation derived from the quadratic form \ref{quad_form_aug} 
\begin{eqnarray}\label{bellman_output}
   \bar{Z}(t)^T\bar{H} \bar{Z}(t)&=& (y(t)-r(t))^TQ(y(t)-r(t)) + u^T(t)Ru(t) \nonumber\\& &   + \gamma \bar{Z}(t+1)^T\bar{H} \bar{Z}(t+1)  .
\end{eqnarray}
To obtain an approximation formula for $\bar{H}$, it needs to be isolated from the quadratic form. This can be achieved by firstly vectorising  \ref{bellman_output}:
\begin{eqnarray}
    vec \left( \bar{Z}^T(t)\bar{H}\bar{Z}(t) \right)& =&  (y(t)-r(t))^TQ(y(t)-r(t)) + u^T(t)Ru(t)\nonumber\\ &&+ \gamma  vec \left( \bar{Z}^T(t+1)\bar{H}\bar{Z}(t+1) \right)
    \end{eqnarray}
    where the $vec$ operator stacks the columns of a matrix to create a single column. We then apply the "vector trick" associated with the Kronecker product
    \begin{eqnarray}\label{pe_step}
    \left( \bar{Z}^T(t)\otimes \bar{Z}^T(t) \right)vec(\bar{H})&= &(y(t)-r(t))^TQ(y(t)-r(t)) +  u^T(t)Ru(t)  \nonumber\\&&+ \gamma  \left(\bar{Z}^T(t+1)\otimes \bar{Z}^T(t+1) \right) vec(\bar{H}) .\label{final_value_func}
\end{eqnarray}
To solve \ref{pe_step}, we may employ the Policy Iteration (PI) or Value Iteration (VI) algorithm \cite{sutton2018reinforcement}. In the standard Reinforcement Learning framework, where the problem is modelled as a Markov Decision Process with discrete state and action spaces, there are two main differences between the algorithms: Firstly, the PI assumes that the initial policy is an admissible one while VI has no such requirement. The second difference lies in the way the optimal policy is built. VI focuses on optimising the value function over all possible actions and extracting the optimal policy at the end. PI updates the value function based on the current policy estimate and extracts a new, better policy afterwards, iteratively until the optimal is reached. In the case of continuous state and action spaces, such as our control problem, we will focus on the VI algorithm, mainly due to the lack of need for an initial stabilising controller (admissible policy). When the system dynamics are not known, the design of such a controller becomes very challenging, though some work has been done in recent literature on obtaining stabilising controllers through measured data \cite{lopez2023efficient}. The algorithm iterates between a Policy Evaluation and a Policy Improvement step. We start with an arbitrarily chosen initial
kernel matrix $\bar{H}^0$ that obtains an initial policy $u^0(t)$. We then iterate between the following two steps: 
\begin{enumerate}
    \item \textbf{Policy Evaluation} 
    \begin{eqnarray*}
        \left( \bar{Z}^T(t)\otimes \bar{Z}^T(t) \right)vec(\bar{H}^{i+1})&=& (y(t)-r(t))^TQ(y(t)-r(t))  +  u^T(t)Ru(t) \nonumber\\&&+ \gamma  \left(\hat{\bar{Z}}^T(t+1)\otimes \hat{\bar{Z}}^T(t+1) \right) vec(\bar{H^i}) 
    \end{eqnarray*}
   
    \item \textbf{Policy Improvement}
    \begin{eqnarray*}
        u^{i+1}(t) &=&- (H_{uu}^{i+1})^{-1} ( H_{u\bar{u}} \bar{u}^{i+1}(t-1,t-N) + H_{u\bar{y}}^{i+1}\bar{y}(t-1,t-N)\\  & & + H_{ur}^{i+1} r(t-N))
     \end{eqnarray*}
\end{enumerate}
The Policy evaluation step can be solved through the Least Squares (LS) algorithm using measured data $\bar{Z}(t)$ and calculating the corresponding cost term $(y(t)-r(t))^TQ(y(t)-r(t)) +  u^T(t)Ru(t)$. It is worth noting that for the discounted future cost term \sloppy $\gamma  \left(\hat{\bar{Z}}^T(t+1)\otimes \hat{\bar{Z}}^T(t+1) \right) vec(\bar{H^i}) $, we are not using the measured input $u(t+1)$ but instead the estimated optimal input $u^i(t+1)$, namely $\hat{\bar{Z}}(t+1)=\begin{bmatrix}
       \bar{u}(t,t-N+1)&  \bar{y}(t,t-N+1) &r(t-N+1) &u^i(t+1)
  \end{bmatrix}^T$. In RL terms, this yields the expected future reward based on the current policy. Matrix $\bar{H}$ is an $((N+1)m+(N+1)p)\times((N+1)m+(N+1)p)$ symmetric matrix which means that its determination is a problem that requires a minimum of $((N+1)m+(N+1)p)\times((N+1)m+(N+1)p)/2$ data points. In practice many more data points are usually needed, especially when dealing with more complicated systems with larger state and action spaces. When the data are highly correlated, which is often the case with sequential input output data from a control system, the inversion step in the LS algorithm becomes challenging. This can be rectified with a regularized LS approach, which includes an appropriate regularisation parameter $\mu$. A step by step description of the method introduced in this section is given in Algorithm \ref{alg_data_driven}.
    
 \begin{algorithm}
		\SetAlgoLined
		\KwResult{Optimal kernel matrix estimate $H^*$, optimal controller estimate $u^*(t)$}
		\KwIn{Weighting matrices for performance index $Q_1,R$, regularisation parameter $\mu$, Reference signal $r(t)$, Discount factor $\gamma$, Initial kernel matrix $H^0$, Error threshold $\epsilon_{rl}$, time horizon $N$}
		\begin{enumerate}
         \item Collect $M \gg ((N+1)m+(N+1)p)\times((N+1)m+(N+1)p)/2$ consecutive data tuples $\left\lbrace y(t),u(t)\right\rbrace$. 
         \item Create vectors $\bar{u}(t-1,t-N), \bar{y}(t-1,t-N)$ and $\bar{Z}(t)=\begin{bmatrix}
       \bar{u}(t-1,t-N)&  \bar{y}(t-1,t-N) &r(t-N) &u(t)
        \end{bmatrix}^T$.
        \item For each $\bar{Z}(t)$, calculate the Kronecker product $K_{\bar{Z}}(t):=\bar{Z}^T(t)\otimes \bar{Z}^T(t)$.
        \item Set $i=0$.
		        \item For t=N to M:
		        \begin{enumerate}[(i)]
		            \item Determine $u^{i}(t+1) =- (H_{uu}^{i+1})^{-1} ( H_{u\bar{u}} \bar{u}^{i+1}(t-1,t-N)+ H_{u\bar{y}}^{i+1}\bar{y}(t-1,t-N) + H_{ur}^{i+1} r(t-N))$
		        \item Create the updated augmented state $\hat{\bar{Z}}(t+1)=\begin{bmatrix}
		            \bar{u}(t,t-N+1)&  \bar{y}(t,t-N+1) &r(t-N+1) &u^i(t+1)
		        \end{bmatrix}^T$ and the corresponding Kronecker product $K_{\hat{\bar{Z}}}(t+1)=\hat{\bar{Z}}^T(t+1)\otimes \hat{\bar{Z}}^T(t+1)$
		        \item Compute the future expected cost term $c(t)=X^T(t)Q_1X(t) $ $ +  u^T(t)R u(t)+ \gamma  K_{\hat{\bar{Z}}}(t+1) vec(H^i)$
		        \end{enumerate}
		        \item Calculate the sums $L=\sum_{t=N+1}^{M} K^T_Z(t)K_Z(t)$ and $R=\sum_{t=N+1}^{M}  K^T_Z(t)c(t)$ .
		        \item Regularise matrix L: $L_r=L+\mu I$
		        \item Produce new H estimate $H^{i+1}=L_r^{-1}R$
		    \item Set $i=i+1$
      \end{enumerate}
      Repeat steps (5)-(9) until $\Vert H^{i}-H^{i+1}\Vert \leq \epsilon_{rl}$ or for a predetermined number of iterations.
		\caption{ \label{alg_data_driven} LQT output-feedback control using input-output measured data}
	\end{algorithm}   

From a software point of view, the implementation of Algorithm \ref{alg_data_driven} can become very computationally expensive with larger system dimensions. There are a few ways to work around this problem. The construction of matrix $L$ in step (6) consists of the Kronecker product terms obtained directly from the dataset. This means that $L$ can be determined only once outside the training iterations. Similarly for its regularised version $L_r$ and more importantly for the inverse $L_r^{-1}$. Matrix inversion is a very expensive operation and an approximation of the inverse might be necessary. Specifically, we chose to calculate the inverse with the use of the Cholesky decomposition \cite{6710599}.

\section{Parameter Tuning through Bayesian Optimisation}\label{section_BO}

There are multiple aspects that influence the performance of a controller, both in the model-based scenario and especially in the model-free, data-driven scenario. In the model-based case, these aspects concern mainly the control objective design. In the case of data-driven control, the resulting controller's performance is influenced as much by characteristics of the data as by the definition of the optimisation problem. In the next section we discuss the setup of our simulations so as to generate informative data, that can be used to effectively train a controller. In this section, we focus on optimising the design of the performance index, which dictates the control objective. 

Apart from the system dynamics \eqref{original_sys}, the main parameters influencing the optimisation are the weighting matrices $Q$ and $R$ in the performance index \eqref{cost_func}. However, there is no systematic way to choose the optimal values for these matrices, as their direct relationship to the controller performance cannot be expressed by a known function. Hence, they are commonly chosen empirically by the engineer, or for simpler problems, exhaustive search may be performed to obtain the best choice. For our system, an exhaustive search would not be feasible, especially in the model-free case where the calculations become vastly more computationally expensive due to the system augmentation. Motivated by this we explore Bayesian Optimisation \cite{frazier2018tutorial,shahriari2015taking,snoek2012practical}, which is a powerful optimization tool in which uncertainty over the objective function is typically represented as a Gaussian process (GP)\cite{rasmussen2006gaussian}. 

In particular, we want to design a fitness function which quantifies the effect of $Q$ and $R$ on the controller's performance. We consider a controller to have good performance if it tracks the reference signal closely, whilst not applying costly inputs to the system. A choice for the fitness function could be the performance index itself as defined in  \ref{cost_func}. There are two main drawbacks to this choice. Firstly, an infinite sum is not feasible, which means it would need to be replaced by a finite sum corresponding to the length of the experiment. Secondly, using \eqref{cost_func} as the fitness function to be minimised, may imply that using very small values for $Q$ and $R$ is sufficient to have a good performance. Such a choice however, would only ensure that the sum of discounted costs is small, without necessarily guaranteeing good tracking performance and low input costs. Hence, we define the simpler fitness function 
\begin{equation}\label{fitness}
    F(\theta)=\sum_{t=0}^l \gamma^t ( \Vert y(t)-r(t)\Vert^2 + \Vert u(t) \Vert^2)
\end{equation}
where $\Vert \cdot \Vert$ denotes the Euclidian norm and $l$ is the length of the experiment. $\theta\in \mathcal{D}$ represents the parameterised choice of $Q$ and $R$ sampled from a bounded search domain $\mathcal{D}$. Considering that $Q$ and $R$ are symmetric positive definite matrices, we can chose them to be diagonal matrices with all positive entries, as is common practice for simplicity. This means that $\theta$ can be written as $\theta= \{ q_1,\dots,q_p,r_1,\dots,r_m  \}$ and $\mathcal{D} \subset \mathbb{R}_+^{p+m}$. This choice for $F$ allows us to investigate the tracking and input costs as an indirect result of $Q$ and $R$, with the weighting matrices being used in the Bellman equation to obtain the optimal policy and thus the optimal trajectory.

We want to minimise the fitness function with respect to $\theta$. $F$ can be modelled through a GP prior as
\begin{equation}\label{prior}
    F(\theta) \sim  \mathcal{GP} (\mu(\theta),k(\theta,\theta'))
\end{equation}
The GP is characterised by a mean $\mu(\theta)$ and a covariance or kernel function $k(\theta,\theta')$. 

Consider a sample of $l$ points $\{ x_1,\dots,x_l \}$ and calculate the covariance for each pair $\{ x_i,x_j \}$. The mean for this initial sample may be chosen as  $\mu_0=0$ which is a common choice for the mean when no prior information is available. For the kernel, we choose one of the most popular and simpler kernel function which is the \textit{Square Exponential} or \textit{Gaussian} kernel
\begin{equation}
    k_0(\theta,\theta')= \exp ( - \Vert \theta - \theta' \Vert^2) \text{ where } \Vert \theta - \theta' \Vert^2 = \sum_{i=1}^{p+m}  (\theta_i-\theta_i')^2.
\end{equation}
All the above yield the prior distribution 
\begin{equation}
    F(x_{1:l})\sim \mathcal{N}(\mu_0(x_{1:l}),k_0(x_{1:l} ))
\end{equation}
where $F(x_{1:l})=\{ F(x_{1}), \dots , F(x_{l})\}$, $ \mu_0(x_{1:l})=0_l$ and $k_0(x_{1:l})=\{ k(x_{1},x_{1}),\dots , $  $ k(x_{1},x_{l}); \dots ; k(x_{l},x_{1}), \dots , k(x_{l},x_{l}) \}$.
Once the prior distribution has been determined for the initial sample of points, we can evaluate $F$ at a new point $x$ and compute the conditional distribution of $F(x)$ using Bayes' rule
\begin{equation}\label{posterior}
    F(x)\vert F(x_{1:l})\sim \mathcal{N}(\mu_l(x),k_l(x ))
\end{equation}
where $\mu_l(x)=k_0(x,x_{1:l})k_0(x_{1:l},x_{1:l})^{-1}(F(x_{1:l})-\mu_0(x_{1:l}))+\mu_0(x)$ and 
$k_l(x)=k_0(x,x)-k_0(x,x_{1:l})k_0(x_{1:l},x_{1:l})^{-1}k_0(x_{1:l},x)$. Distribution \eqref{posterior} is also known as the \textit{posterior probability distribution}. The more points that are sampled, the more informed the posterior becomes. 

Returning to the optimisation problem at hand, we are looking for a choice of $Q$ and $R$ that minimise $F$. The tool used in Bayesian Optimisation to find such a point is known as the \textit{acquisition function}. Acquisition functions provide a candidate for the best point to sample next, based on the current estimate for the distribution. There are a few different choices for acquisition functions, the most common including \textit{Expected Improvement, Probability of Improvement} and \textit{Gaussian Process Upper Confidence Bound}. For our application, we will be using the latter, defined as follows
\begin{equation}
    a_{UCB}(x;x_{1:l}) =\mu_l(x)-\epsilon k_l(x) 
\end{equation}
where $\epsilon$ is a parameter balancing exploration against exploitation. The next point to sample $V$ is chosen as $x_{next}=\arg\min a_{UCB}$. The new point can be used to update the posterior distribution, while resulting in a more informed acquisition function that will in turn suggest a better point to sample next. The process is repeated for a predefined number of steps and the optimal is extracted at the end if $M$ points have been sampled in total, as 
\begin{equation}
    x_{best}=\arg\min_{x}\mu_M(x).
\end{equation}

\section{Simulation setup and results} \label{section_results}
In this section we perform numerical simulations to investigate the behaviour of the extruder heating system introduced in Section \ref{section_sys_intro} after a controller is introduced to the system. We assume partial access to the system state in the form of measurements and explore both the case of full knowledge of the system dynamics and the case where a system model is unknown. For our simulations, we need a model of the extruder dynamics which will be used for the design of the model-based controller and to generate data to train the data-driven controller. In \cite{gootjes2017applying}, system identification methods were used to obtain the following state and input matrices for a BAAM system:

\begin{eqnarray} \footnotesize \label{sys_mats_ab}
    A&=& \begin{bmatrix}
    0.992 & 0.0018 & 0 & 0 & 0 & 0 \\ 
    0.0023 & 0.9919 & 0.0043 & 0 & 0 & 0 \\ 
    0 & -0.0042 & 1.0009 & 0.0024 & 0 & 0 \\ 
    0 & 0 & 0.0013 & 0.9979 &0&0 \\ 
    0&0&0&0&0.9972&0 \\
    0&0&0&0&0&0.9953
    \end{bmatrix} \nonumber\\ 
    B&=& \begin{bmatrix}
    1.0033 & 0 & 0 & 0 & 0 & 0&-0.2175 \\ 
    0 & 1.0460 & 0 & 0 & 0 & 0&-0.0788 \\ 
    0 & 0 & 1.0326 & 0 & 0 & 0 &-0.0020\\ 
    0 & 0 & 0 & 0.4798 &0&0 &-0.0669\\ 
    0&0&0&0&0.8882&0 &0.1273\\
    0&0&0&0&0&1.1699&-0.1792
    \end{bmatrix} \nonumber \\
\end{eqnarray}
In addition, we design an output matrix for the system as follows:
\begin{equation} \label{sys_mats_c}
    C = \begin{bmatrix}
    0.992&0.00018&0&0&-0.0001&0\\
              0.0023&1.3&0.0043&0&0&0\\
              0&-0.0042&1.0109&0.0024&0&0.201\\
              0&0&0.0013&0.989&0.00031&0.64\\
              0&0&0&0&0.923&0.3
    \end{bmatrix}
\end{equation}
We can easily verify that the system described by the above matrices is controllable and observable. As highlighted in the previous section, the most important design parameters for the optimisation are the weighting matrices Q and R. Before performing Bayesian Optimisation to obtain the best values, we chose as a baseline identity matrices of appropriate dimensions, namely $Q=I_p$ and $R=I_m$.

In BAAM systems, the multiple heaters located across the barrel heat the material in steadily increasing temperatures until the material reaches a suitable temperature for extrusion. For our simulations, we assume the polymer used is PLA, whose melting temperature ranges between 160 and 180 \textdegree C. With this in mind, we set the reference signal to be $r(t) = \{  150,160,170,175,180 \}, \forall t\geq t_0$. This means that the reference generating matrix is $F=I_5 $. This reference will be used throughout the simulations, both for the model-based controller and the data-driven one. 

\subsection{State Observer based Output Feedback}\label{section_results_observer}

We first focus on the approach introduced in Section \ref{section_model_based}, starting with the design of the State Observer, for which we choose $\tau=0.002$. For the controller design, we use a discount factor $\gamma=0.99$ and an error threshold for convergence $\epsilon=0.01$. We arbitrarily initialise the state estimate with $\hat{x}(t_0)= \{ 50,\dots, 50 \}$ and the control policy with random normally distributed numbers as $K^0=\{ k_1, \dots , k_{n+p}\}$ where $k_i\sim \mathcal{N}(0,10)$. Finally, even though we assume that the system state is not directly available, we need to initialise it so that it can be updated and measured during the simulations. That starting point is chosen to be $x(t_0)= \{ 20,\dots, 20 \}$. Table \ref{table_observer_based} summarises all parameter values used for the observer-based simulations.
	\begin{table}[ht]
	\centering
	\begin{tabular}{|c|c||c|} 
	\hline
	\multicolumn{3}{|c|}{Table of parameters - Observer based output feedback}\\
		\hline
		\hline
		Parameter&Symbol&Value\\
		\hline 
		Tracking reference & r & $\left\lbrace 150,160,170,175,180 \right\rbrace$ \\
        Initial state & $x(t_0)$ & $\left\lbrace 20, \dots, 20 \right\rbrace$\\
		Initial observer state & $\hat{x}(t_0)$ & $\left\lbrace 50, \dots, 50 \right\rbrace$\\
        Observer pole & $\tau$ & $0.002$ \\
        Error threshold &$\epsilon$ & $0.01$\\
        Initial control gain & $K^0=\{ k_i\}_{i=1}^{n+p}$&  $k_i\sim \mathcal{N}(0,10)$ \\ 
		\hline
	\end{tabular} 
	\caption{\label{table_observer_based}Parameters chosen for observer-based control implementations. The information listed includes the parameter names, the corresponding symbol and the value that was chosen. }
\end{table}

\begin{figure}[ht]
    \centering
    \includegraphics[width=0.7\linewidth,height=5cm]{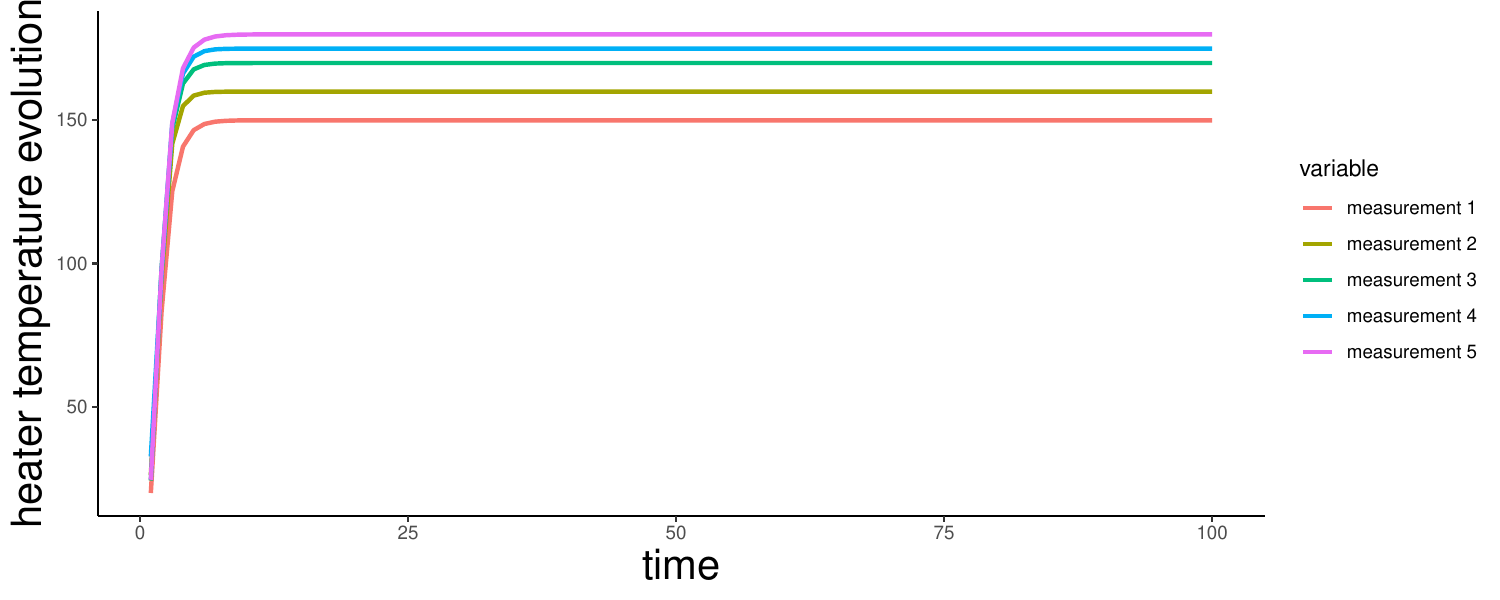}
    \caption{Observer-based Output Feedback. The system output achieves temperatures very close to the optimal within a few steps. }
    \label{fig:observer_trajectories}
\end{figure}

\begin{figure}[ht]
	\begin{subfigure}{.45\textwidth}
		\includegraphics[width=\linewidth,height=4cm]{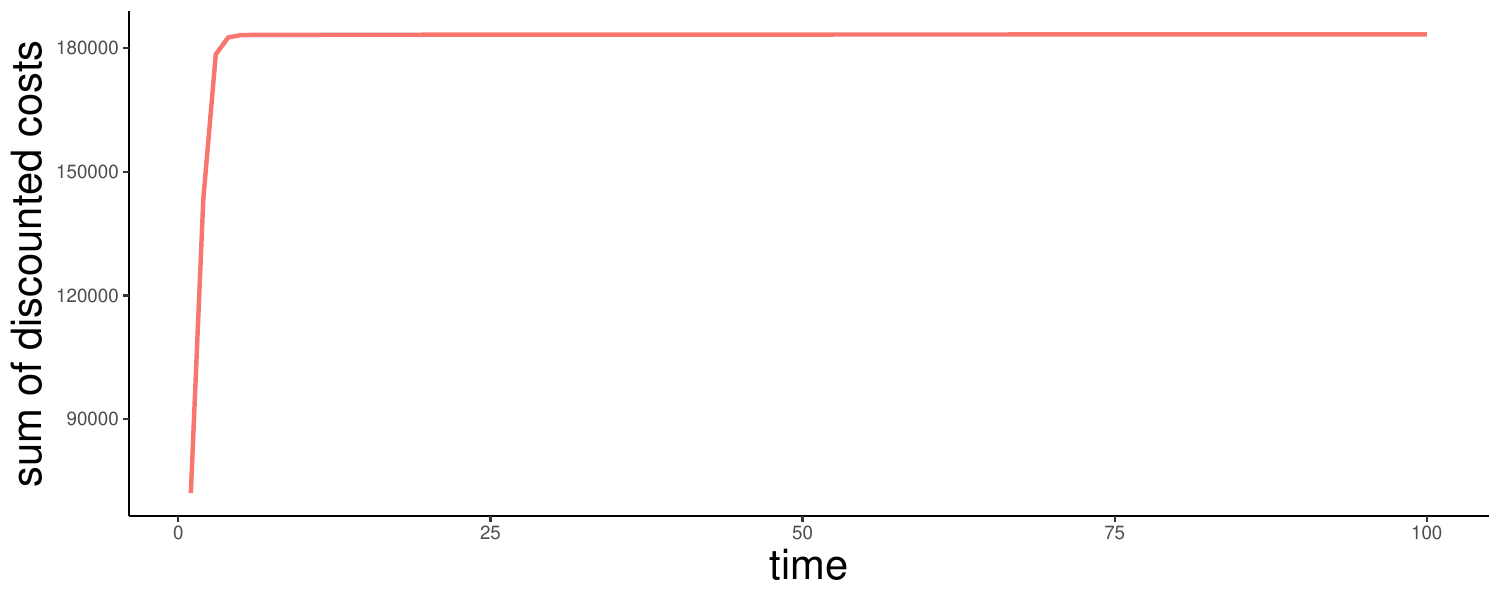}
		\caption{\label{fig:observer cost}Sum of discounted cost terms for 100 time steps.}
	\end{subfigure}%
 \hspace*{\fill}
	\begin{subfigure}{.45\textwidth}
		\includegraphics[width=\linewidth,height=4cm]{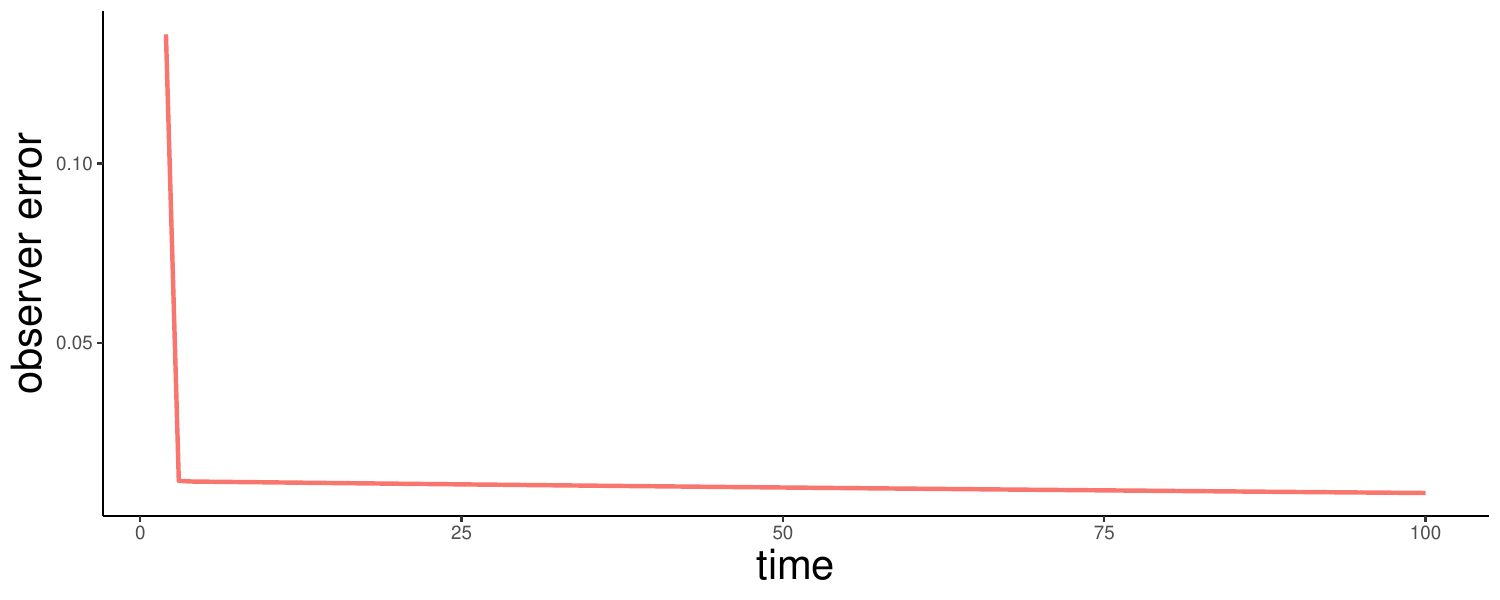}
		\caption{\label{fig:observer error}Estimation error for the system state by the observer.}
	\end{subfigure}
	\caption{\label{fig:observer cost and error}Performance index and state estimation error for the Observer-based Output feedback. }
	\end{figure}

Figure \ref{fig:observer_trajectories} shows the output trajectories obtained as a result of applying the model-based controller for 100 time steps, as described in Algorithm \ref{observer_algorithm}. The controller manages to bring all measurements within $0.03$\textdegree C of the optimal. The measured temperatures arrive within $0.1$ degrees of the target within only 10 steps and remain there. The values that the measurements finally converge to are $y^*=\{   149.9798,      159.9932,      169.9795,      174.9815  ,    179.9859\}$. The error between the final measurement values and the reference, calculated through the Euclidean norm, is $\Vert r-y^*\Vert_2=0.0376$. The value of the performance index obtained by this controller, approximated as the discounted sum of the first 100 cost terms is 183,362.5. When calculated over 1000 steps, it reaches 183.436.2. The performance index trajectory for 100 time steps can be seen in Figure \ref{fig:observer cost}. The good performance of the controller is due to the fact that the observer approximates the system state very closely, thanks to the knowledge of the system dynamics. The estimation error, calculated as the Euclidean norm $e(t)=\Vert\hat{x}(t)-x(t)\Vert_2$, is decreasing with each time step, reaching the value $0.008$ in 100 steps. It is worth noting that the closer the initial state estimate is to the initial state, the faster the error decreases. This decreasing trend can be seen in Figure \ref{fig:observer error} where we plot the estimation error for the first 100 time steps. 

We perform Bayesian Optimisation to obtain better choices for the weighting matrices $Q$ and $R$ and study the effect they have on our results. For simplicity, as is common practice, we assume these matrices are diagonal and positive definite. Through some trial and error attempts, we find that the best range of values for the diagonal entries is positive numbers $\leq 1$. To avoid numerical instability, we assume a lower limit for these values to be $0.01$. Additionally, we want to put a higher emphasis on the tracking of the reference signal. This is expressed by choosing appropriate upper limits for the entries of each matrix. Specifically, we chose $1$ as the upper limit for $Q$ and $0.4$ as the upper limit for $R$. 

\begin{figure}[ht]
	\begin{subfigure}{.45\textwidth}
		\includegraphics[width=\linewidth,height=4cm]{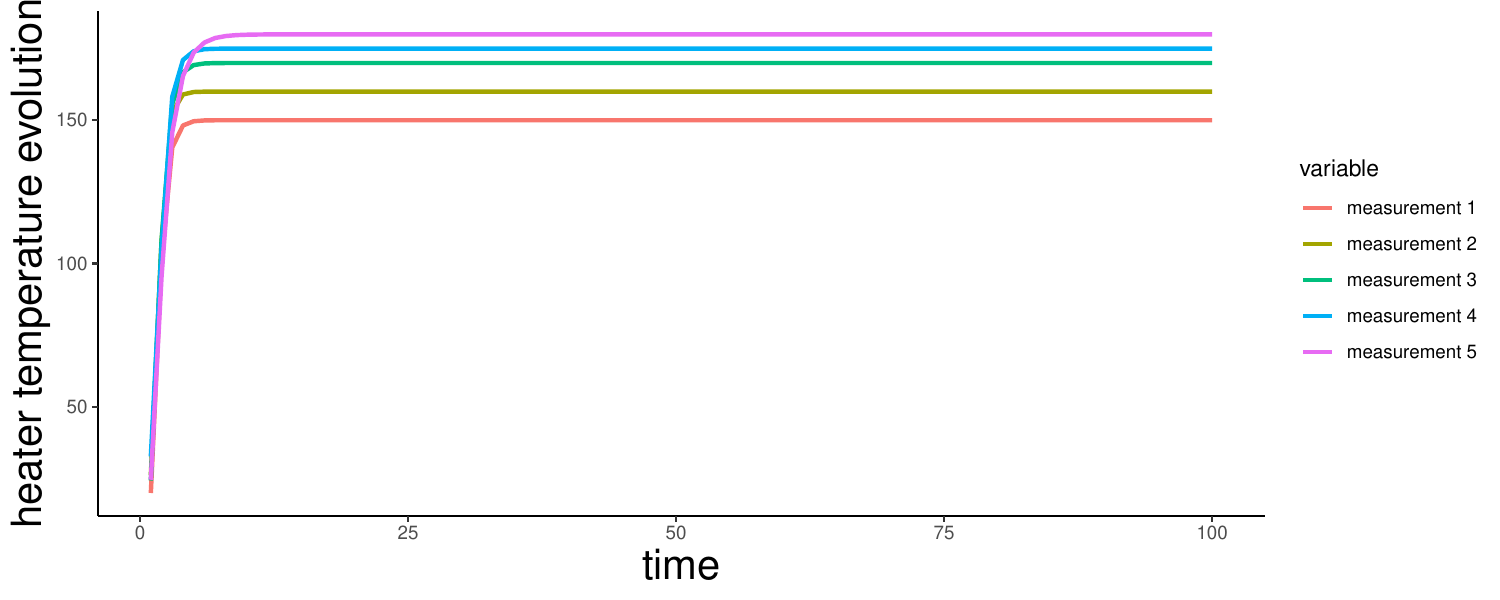}
		\caption{\label{fig:observer BO trajectories}System output trajectories}
	\end{subfigure}%
 \hspace*{\fill}
	\begin{subfigure}{.45\textwidth}
		\includegraphics[width=\linewidth,height=4cm]{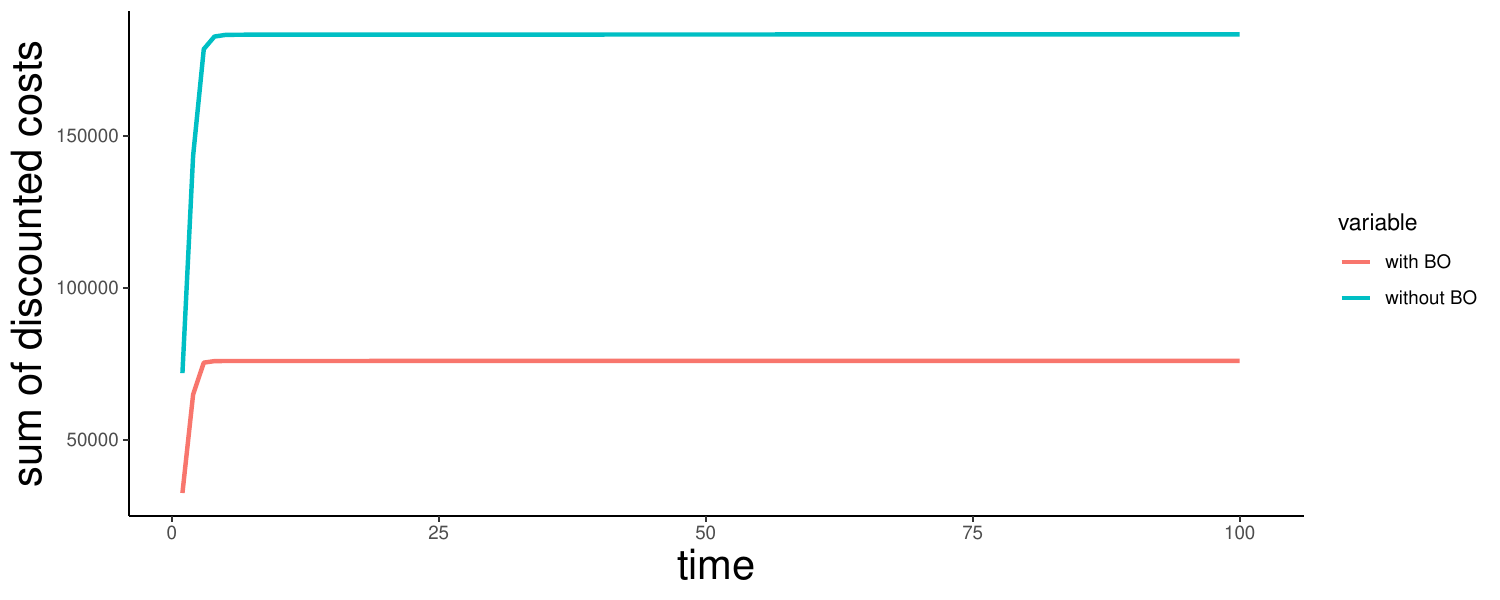}
		\caption{\label{fig:observer BO error}Comparison of performance index before and after using Bayesian Optimisation}
	\end{subfigure}
	\caption{\label{fig:observer BO}Output trajectories and performance index comparison for the Observer-based Output feedback where the weighting matrices $Q$ and $R$ were obtained through Bayesian Optimisation. The output trajectories are similar to the default case while the performance index is minimised more effectively through the Bayesian Optimisation approach.}
	\end{figure}

We run the Bayesian Optimisation procedure as described in Section \ref{section_BO} with 50 initial randomly sampled points for 100 iterations. The optimal values obtained for the matrices Q and R rounded up to 3 decimal points are 
\begin{equation}
    Q^*=\begin{bmatrix}
        0.943 & 0 & 0& 0 & 0 \\
        0 & 0.762& 0& 0&0 \\
        0 & 0 & 0.542 & 0 & 0 \\
        0 & 0 & 0 & 0.420 & 0 \\ 
        0 & 0 & 0 & 0 & 0.514
    \end{bmatrix}
\end{equation}
and 
\begin{equation}
    R^* = \begin{bmatrix}
        0.300 & 0 & 0& 0 & 0 & 0 & 0 \\ 
        0 & 0.270 & 0 & 0 & 0 & 0 & 0 \\
        0 & 0 & 0.281 & 0 & 0 & 0 & 0 \\
        0 & 0 & 0 & 0.092 & 0 & 0 & 0 \\
        0 & 0 & 0 & 0 & 0.054 & 0 & 0 \\ 
        0 & 0 & 0 & 0 & 0 & 0.269 & 0 \\
        0 & 0 & 0 & 0 & 0 & 0 &  0.318
    \end{bmatrix}
\end{equation} 
Using the optimised values $Q^* $ and $R^*$ we obtain the updated trajectories seen in Figure \ref{fig:observer BO trajectories}. The values that the outputs converge to are: $\{  149.9940 ,     159.9946 ,     169.9843 ,     174.9873   ,   179.9834\}$ which means they all arrive within $0.02$\textdegree C of the optimal (cf $0.03$\textdegree C for the previous version). The error of the final values from the optimal is $\Vert r-y^*\Vert_2=0.0274$, verifying that the tracking performance is improved. Additionally, the performance index is significantly lower in this case, as can be seen in Figure \ref{fig:observer BO error}. Specifically, the performance index obtains the value 76,060.36 after 100 steps and converges to 76,072.6 over 1000 steps. This means that the updated weighting matrices result in a $58.5 \%$ reduction of the performance index.

\subsection{Data generation} \label{section_data_generation}
When designing and training model-free, data-driven controllers, much like any learning algorithm, it is important to have data that is rich in information. In the case of the algorithm introduced in Section \ref{section_data_driven}, the iterative algorithm learns the kernel matrix $H$ through input-output tuples. Matrix $H$ includes information of the system dynamics as seen in Equation \ref{kernel}, meaning that learning $H$ is equivalent to learning the system model. This means that the problem can also be viewed as a system identification problem \cite{ljung1998system}. 

An important advantage of working with simulated data generated from an assumed model is that we can ensure the data is appropriate to train our algorithms. The way to achieve that is by using a well designed input signal to the system. Firstly, we need the system output to stay bounded and not diverge towards infinite values, in which case the data would not be usable. To that end, we can design the input to be in feedback form \begin{equation}\label{feedback_input}
    u(t)=-K_{data} x(t)
\end{equation} 
where the feedback gain $K_{data}$ is stabilising. In the case where the training data is obtained from a real process, this step is not applicable, however, the process would be tuned from the manufacturer to avoid extreme divergence. In the case of the system introduced in \ref{sys_mats_ab} and \ref{sys_mats_c}, we can use the \textit{Matlab Control System Toolbox} to obtain a stabilising gain. In particular, two of the most popular approaches to stabilising a system's behaviour are the Pole Placement method and the Linear Quadratic Regulator (LQR). Using the LQR approach, we obtain the gain:
\begin{equation}\label{data_sim_K}
K_{data}= \begin{bmatrix*}[r]
    0.7395& -0.0076& -0.0003&  -0.0264&    0.0194& -0.0170\\
    -0.0076&   0.7430 &   0.0031&   -0.0093&  0.0068&  -0.0060\\ 
    -0.0003&   -0.0033&    0.7599&    0.0021&    0.0002&  -0.0002\\  
    -0.0126&   -0.0042&   0.0016&   1.0971&   0.0092&  -0.0079\\ 
    0.0171&    0.0058&   0.0002&   0.0170&    0.8179&    0.0108\\ 
    -0.0198&   -0.0067&  -0.0002&  -0.0193&    0.0143&    0.6823\\ 
    -0.1525&   -0.0519&   -0.0018&   -0.1412&    0.1091&  -0.0977
\end{bmatrix*}  
\end{equation}
Secondly, in order to ensure the output data is varied enough to efficiently train the control algorithm, we require the input signal to not be `too simple'. This requirement can be fulfilled by the addition of an appropriate probing noise signal to the input, meaning that it will now be of the form
\begin{equation}\label{feedback_with_noise}
    u(t)=-K_{data} x(t) + \omega_{pr}. 
\end{equation} 
We define the probing noise as a sum of a Gaussian noise signal and sinusoidal signals, of varying ranges and frequencies. Specifically, to ensure great variety in the data, we chose ranges between 10 and 100. The choice of sampling frequencies is more particular, where an appropriate choice is the system's bandwidth, i.e. the part of frequency domain where most of signal's energy is contained. Formally, a system's bandwidth is the frequency at which the magnitude of the system's response drops below 0.707 (-3dB). In the case of multi-input multi-output (MIMO) systems, the bandwidth is the frequency where the maximum singular value of the frequency response crosses 0.707 from below \cite{shinners1998modern}. In Matlab, this value can be easily obtained for any system using the \textit{sigma} function. The bandwidth for our system is 1.65, so we use it as the maximum frequency from which to sample sinusoidal signals. Finally, we add a normally distributed random vector to our noise signal, to add an extra element of randomness. We consider the variance for this vector to be $ \sigma =1.5$. Bringing all the above elements together, we have
\begin{eqnarray}\label{noise}
    \omega_{pr}&=&\omega_1 +100 \sin(\omega_2 t) +90 \sin(\omega_3 t)+  80 \sin(\omega_4 t) +70 \sin(\omega_5 t)\nonumber \\ 
    && +60 \sin(\omega_6 t)+50 \sin(\omega_7 t)+40 \sin(\omega_8 t)+30 \sin(\omega_9 t) \nonumber \\ &&+20 \sin(\omega_{10} t)+10 \sin(\omega_{11} t)
\end{eqnarray}
where 
\begin{equation}
    \omega_1\sim \mathcal{N}(0,\sqrt{\sigma}I_7), \; \omega_2=(16.5,\dots,16.5)\in \mathbb{R}^7, \text{ and } \omega_3,\dots \omega_{11} \sim \mathcal{U}_7(0,1.65). 
\end{equation}

\subsection{Data-driven Output Feedback }
Using the input signal designed in the previous section, we can obtain data tuples $\{y(t),u(t)\}$ to train the data-driven RL based controller introduced in Section \ref{section_data_driven}. In particular, we generate $M=15000$ samples and set a maximum number of iterations to be 1000 in case the algorithm does not converge according to the error threshold chosen as $\epsilon_{rl}=0.001$. As the system is fully observable, we choose the time horizon $N$ for the augmented data tuples to be $N=6$. We choose the discount factor to be $\gamma=0.99$ and the regularisation parameter to be $\mu=0.0001$. Finally we initialise the kernel matrix $H$ as the identity matrix of appropriate dimensions, namely $H^0=I_{(N+1)m+(N+1)p}$. We run two sets of simulations, one with the weighting matrices $Q, R$ chosen to be identity matrices of appropriate dimensions, and one using the values for $Q $ and $R$ obtained through Bayesian Optimisation. 

	\begin{table}[ht]
	\centering
	\begin{tabular}{|c|c||c|} 
	\hline
	\multicolumn{3}{|c|}{Table of parameters - Model based methods}\\
		\hline
		\hline
		Parameter&Symbol&Value\\
		\hline
		\multicolumn{3}{|c|}{Universal parameters}\\
		\hline 
		Tracking reference & r & $\left\lbrace 150,160,170,175,180 \right\rbrace$ \\
        Initial state & $x(t_0)$ & $\left\lbrace 50, \dots, 50 \right\rbrace$\\
        Discount factor &$\gamma$ & $0.99$ \\ 
        Error threshold & $\epsilon_{rl}$ & $0.001$ \\
        Regularisation parameter & $\mu$ & $0.0001$\\
        Time horizon & $N$ & $6$ \\
        Initial kernel matrix & $H^0$ & $I_{(N+1)m+(N+1)p}$\\
		\hline
	\end{tabular} 
	\caption{\label{table_rl}Parameters chosen for data-driven RL based implementations. The information listed includes the parameter names, the corresponding symbol and the value that was chosen.}
\end{table}

\begin{figure}[ht]
	\begin{subfigure}{.45\textwidth}
		\includegraphics[width=\linewidth,height=4cm]{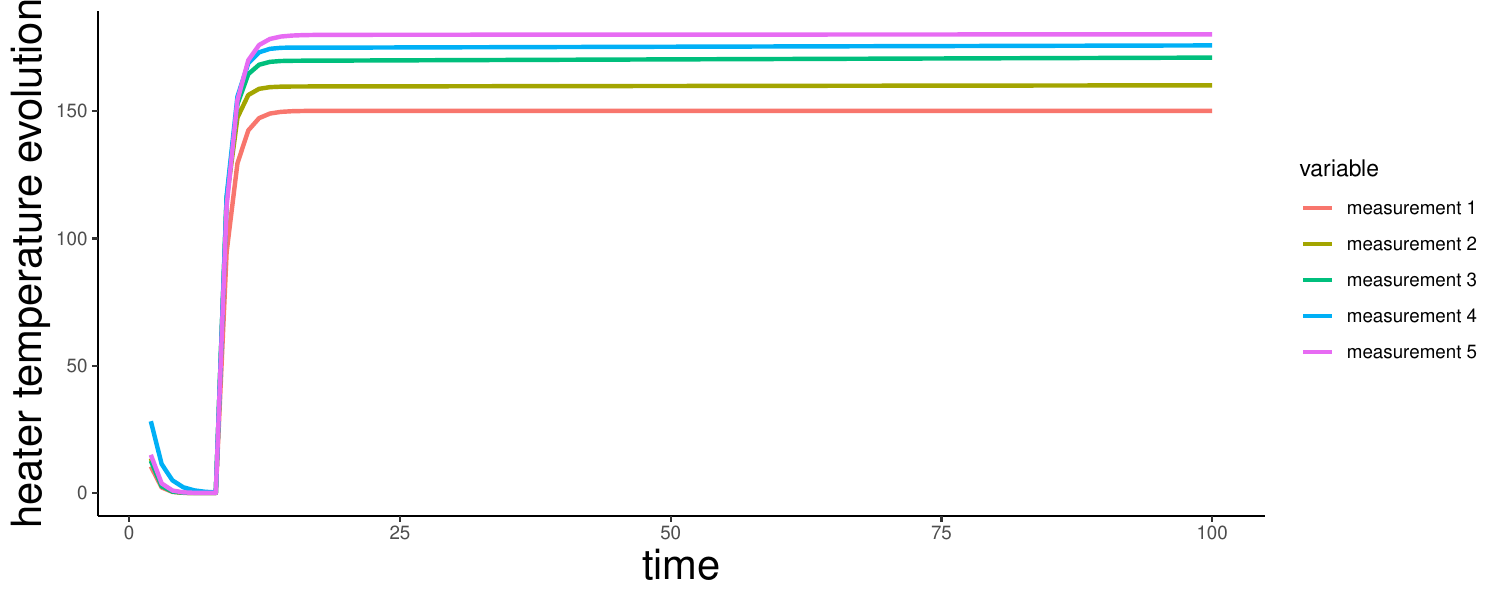}
		\caption{\label{fig:output feedback trajectories no bo}System output trajectories}
	\end{subfigure}%
 \hspace*{\fill}
	\begin{subfigure}{.45\textwidth}
		\includegraphics[width=\linewidth,height=4cm]{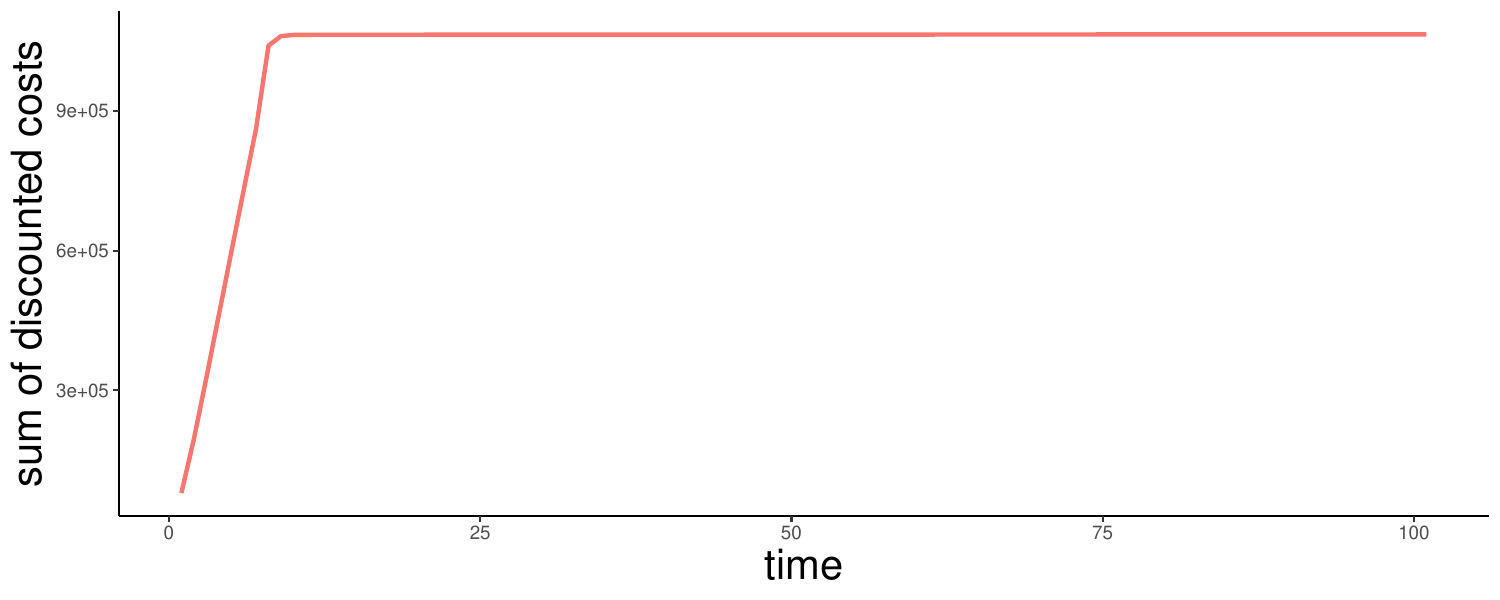}
		\caption{\label{fig:perf index no bo}Performance Index}
	\end{subfigure}
	\caption{\label{fig:output feedback no BO}Output trajectories and the performance index the data-driven Output feedback algorithm using input-output data, where the weighting matrices $Q$ and $R$ were chosen as identity matrices of appropriate dimensions.  }
	\end{figure}

Firstly, we focus on the case of the arbitrarily chosen weighting matrices $Q$ and $R$. We chose them to be identity matrices of appropriate dimensions, namely $Q=I_5$ and $R=I_7$. Figure \ref{fig:output feedback no BO} shows the output trajectories obtained through the data-driven controller described in Algorithm \ref{alg_data_driven}, and the corresponding performance index. The measurements converge to values within $1$\textdegree C of the optimal, specifically to the vector $y^*=\{ 150.0857 , 160.1181,  171.0027,   175.7926,  180.1416\}$. The error between these values and the optimal ones, calculated through the Euclidean norm, is $\Vert r-y^*\Vert_2=1.2942$. The performance index, defined as the infinite sum of discounted costs, is estimated by a large finite sum which for 1000 steps is 1,074,985. We plot it for the first 100 steps, in which period the sum has obtained the value 1,065,077. 

 We now perform Bayesian Optimisation on $Q$ and $R$ using the same upper and lower bounds for the search space as in Section \ref{section_results_observer}, namely we chose $0.1$ as the lower bound for both matrices and chose $1$ as the upper bound for $Q$ and $0.4$ as the upper bound for $R$. We generate a grid of 50 random initial samples and then train the algorithm for 100 iterations. The optimal values found for $Q$ and $R$ rounded to 3 decimal points are 
 \begin{equation}
    Q^*=\begin{bmatrix}
        0.174 & 0 & 0& 0 & 0 \\
        0 & 0.056 & 0& 0&0 \\
        0 & 0 & 0.010 & 0 & 0 \\
        0 & 0 & 0 & 0.010 & 0 \\ 
        0 & 0 & 0 & 0 & 0.160
    \end{bmatrix}
\end{equation}
and
\begin{equation}
    R^* = \begin{bmatrix}
        0.145 & 0 & 0& 0 & 0 & 0 & 0 \\ 
        0 & 0.389 & 0 & 0 & 0 & 0 & 0 \\
        0 & 0 & 0.020 & 0 & 0 & 0 & 0 \\
        0 & 0 & 0 & 0.116 & 0 & 0 & 0 \\
        0 & 0 & 0 & 0 & 0.099 & 0 & 0 \\ 
        0 & 0 & 0 & 0 & 0 & 0.316 & 0 \\
        0 & 0 & 0 & 0 & 0 & 0 &  0.010
    \end{bmatrix}.
\end{equation} 

\begin{figure}[ht]
	\begin{subfigure}{.45\textwidth}
		\includegraphics[width=\linewidth,height=4cm]{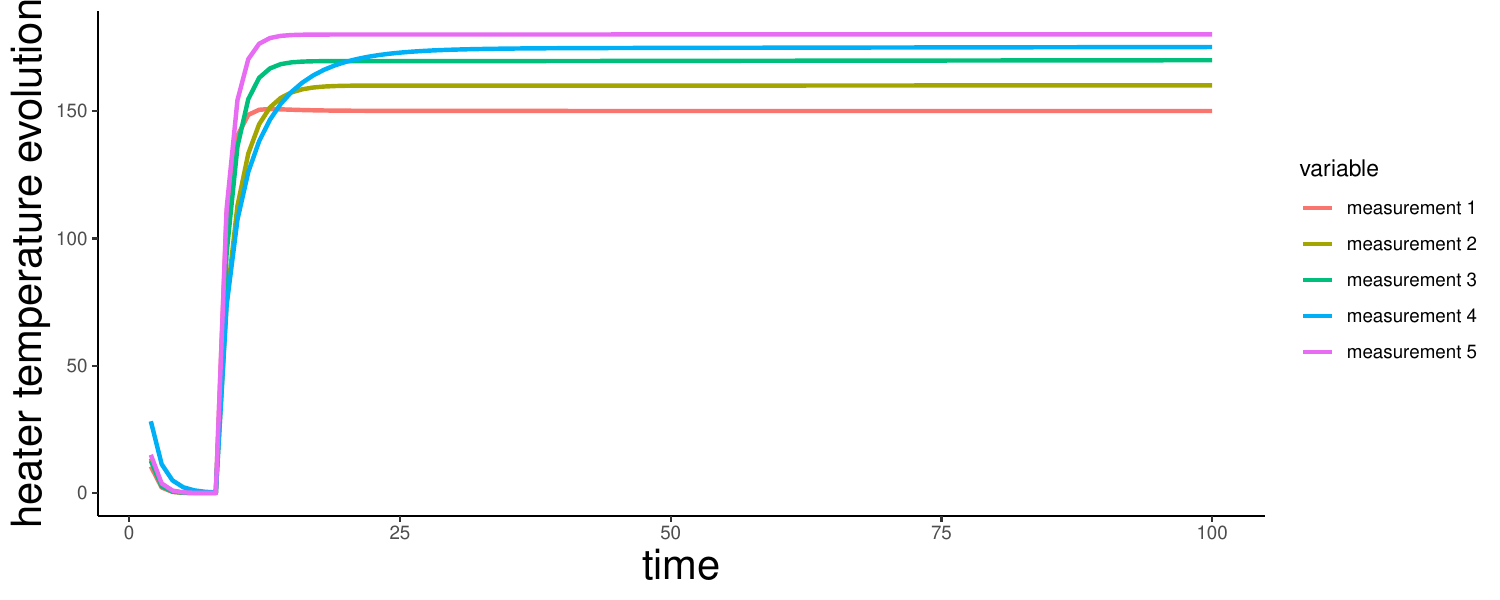}
		\caption{\label{fig:output feedback trajectories bo}System output trajectories after Bayesian Optimisation}
	\end{subfigure}%
 \hspace*{\fill}
	\begin{subfigure}{.45\textwidth}
		\includegraphics[width=\linewidth,height=4cm]{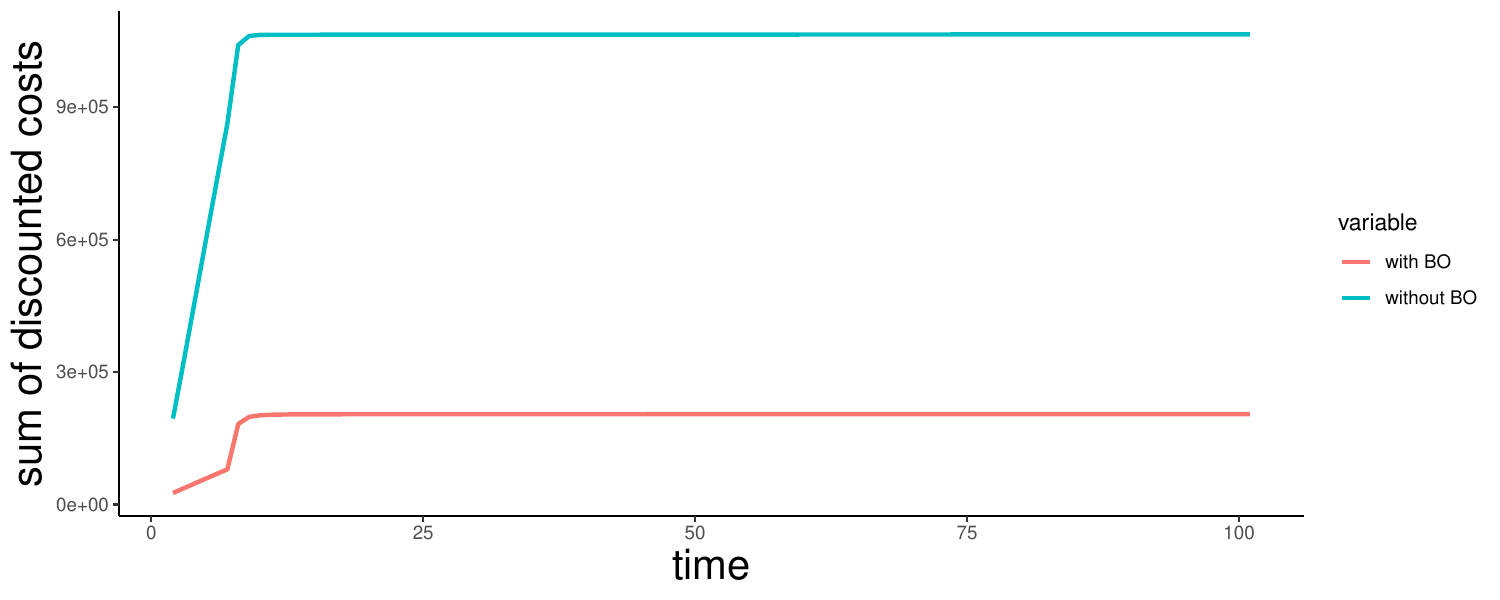}
		\caption{\label{fig:perf index comp output feedback}Cost Comparison before and after Bayesian Optimisation}
	\end{subfigure}
	\caption{\label{fig:output feedback BO} (a) Output trajectories where the weighting matrices $Q$ and $R$ were obtained through Bayesian Optimisation. (b) A comparison of the performance index before and after the optimisation of $Q$ and $R$. There is a significance decrease of the total cost after Bayesian Optimisation, while simultaneously .}
	\end{figure}

With these values for $Q$ and $R$, Algorithm \ref{alg_data_driven} obtains a controller that produces the output seen in Figure \ref{fig:output feedback trajectories bo}. The measurements converge to the vector $y^*= \{ 150.069 160.154  170.011   175.199  180.176\}$ after 100 steps, namely they converge within $0.2$\textdegree C of the optimal. The error associated with these values is $\Vert r-y^* \Vert = 0.3154$. This means that a tracking performance is significantly improved compared to the previous case. In addition, the overall cost is much lower, as can be seen by the comparison of the two performance indices, before and after Bayesian Optimisation, in Figure \ref{fig:perf index comp output feedback}. In particular, the new infinite sum of discounted costs, when approximated by the first 100 steps, obtains the value $204,635.7$, while the value reached after 1000 steps is $206,358.8$. This means that Bayesian Optimisation achieves an 80.8\% decrease in the cost, as expressed through the performance index.

Different implementations might yield slightly different numerical results, depending on the randomness introduced into the algorithm. In particular, randomness is introduced in the data generation step, as well as the sampling for Bayesian Optimisation. However, as long as the learning algorithm is designed according to Algorithm \ref{alg_data_driven} and the data is generated appropriately, as described in Section \ref{section_data_generation}, the resulting controller will successfully track the reference signal, and Bayesian Optimisation will consistently reduce the cost while maintaining good tracking performance. 

\section{Conclusion}\label{section_conclusion}
In this work, we investigated the problem of temperature LQT control within the extruder of a Big Area Additive Manufacturing system, through tracking a reference signal. In particular, we focused on the case where direct access to the system state is not possible, instead some measurements of the state are available. We considered a state-space model identified and validated based on real data describing the system dynamics, and used it to simulate the system's behaviour and response to different controllers. The lack of access to the system state introduces the need for a state estimation step, to better inform the controller. Assuming the system model is available, the task of state estimation can be performed through the implementation of a state observer. The resulting state estimate is then used to inform a feedback controller, designed using model-based Reinforcement Learning Techniques, with the Bellman equation being the main tool. 

Next, we investigated solving the Output Feedback problem in a model-free and data-driven way. We used the state space temperature model for data generation but no longer assumed it available for the control design. We instead used Least Squares to train a controller directly from process data, while incorporating a state estimation procedure in the algorithm. As with many learning algorithms, the data quality is really important when training data-driven controllers. For this reason, we provide an in depth discussion on the input signals that would generate "good" data when applied to the system. Finally, we explored the use of Bayesian Optimisation to improve the choice of weighting parameters in the performance index. 

We presented our results before and after the implementation of Bayesian Optimisation, both for the model-based and data-driven controllers. For the model-based scenario, the obtained controller successfully tracks the reference signal, which was expected given the access to full information about the system. We were able to achieve comparable results while alleviating the need for knowledge of the system dynamics. This is a very encouraging result, considering how challenging it is to have accurate models of AM processes. Our results also indicate the improvement Bayesian Optimisation offers, as expressed through the significant reduction of the performance index, namely a 58.5\% improvement in the model-based case and an 80.8\% improvement in the data-driven case. 

While our approach is a valuable tool to solving the LQT problem, it presents certain limitations. The basis of our approach is the quadratic form of the objective function (performance index). When defining a control problem, the objective function can take many forms, especially when introducing constraints into the problem. For example, we can impose constraints on the output response of the system so that the values follow steadily increasing trajectories and the dip seen in Figures \ref{fig:output feedback trajectories bo} and \ref{fig:output feedback trajectories no bo} is avoided. Another important drawback of our approach stems from its mathematical and computational complexity. In particular, the implementation of Least Squares for the estimation of kernel matrix $H$ involves matrix operations between very large matrices. Normally, the more complex the system to be controlled is, the larger the state space needed to describe it is. This means that the matrices involved in the introduced method can become significantly larger than the ones in our case study.

The above concerns can be addressed through the introduction of more abstract Reinforcement Learning algorithms, based on incorporating Deep Neural Networks \cite{bertsekas1996neuro}. These are powerful function approximators that can handle large amounts of data, thus being able to handle much more complicated systems and alieviating the need for a particular form to the value function. These approaches are widely known as Deep Reinforcement Learning methods. A very popular one among them is the Actor-Critic framework \cite{konda1999actor}, where one network approximates the value function (critic) and another produces the new action to be taken (actor). We leave these approaches as future work which we explore elsewhere.

  \bibliographystyle{elsarticle-num} 
\bibliography{bibl}

\end{document}